\documentclass[10pt,twoside]{article}
\usepackage{appendix}
\usepackage{graphicx}
\usepackage{amsmath}
\usepackage{amsthm}
\usepackage[byname]{smartref}
\usepackage{txfonts}
\usepackage{tocloft}
\usepackage{setspace}
\usepackage{geometry}
\usepackage{amssymb}
\usepackage{cancel}
\usepackage{mathrsfs}
\usepackage{upgreek}
\usepackage[usenames,dvipsnames]{xcolor}
\usepackage{tikz}
\usepackage[english]{babel}
\usepackage{textpos}
\usepackage{hyperref}
\hypersetup{
    colorlinks,
    citecolor=black,
    filecolor=black,
    linkcolor=black,
    urlcolor=black
}
\makeatletter
\addtoreflist{section}

\newtheorem{thm}{Theorem}[section]
\newtheorem{lem}[thm]{Lemma}
\newtheorem{prop}[thm]{Proposition}
\newtheorem{cor}[thm]{Corollary}

\newtheorem{defn}[thm]{Definition}
\newtheorem{eg}[thm]{Example}

\newtheorem*{con*}{Conjecture}
\newtheorem*{ques*}{Question}
\providecommand{\abs}[1]{\lvert#1\rvert}
\providecommand{\Ren}[0]{\mathcal{R}}
\providecommand{\GJ}[0]{\mathcal{GJ}}
\providecommand{\JG}[0]{\mathcal{JG}}
\providecommand{\Oh}[0]{\mathcal{O}}
\providecommand{\NU}[0]{\mathcal{N}}
\providecommand{\Tc}[0]{\mathscr{T}}
\providecommand{\Uc}[0]{\mathscr{U}}
\providecommand{\Jc}[0]{\mathscr{J}}
\providecommand{\Lc}[0]{\mathscr{L}}
\providecommand{\Rc}[0]{\mathscr{R}}
\providecommand{\Hc}[0]{\mathscr{H}}
\providecommand{\Dc}[0]{\mathscr{D}}

\providecommand{\remove}[0]{\backslash}

\begin{document}

\title{The Adherence Order on Renner-Coxeter Monoids}

\author{Allen O'Hara}
\date{January 2016}
\maketitle

\begin{abstract}
Building upon the previous Renner-Coxeter system of work by Eddy Godelle we introduce the familiar Renner monoid structure of the Adherence order. The Green's relations of the system are then considered in relation to the Adherence order and in the finite case maximum and minimum elements in each equivalence class are distinguished and studied.
\end{abstract}

\vspace{5pt}

\section{Coxeter Systems and Double Cosets}

Renner monoids are a much studied combinatorial object coming from reductive algebraic monoids. Named after Lex Renner, they are derived from reductive algebraic monoids by means of the Bruhat decomposition (\cite{Renner Bruhat I}). They occupy a position in the theory of reductive monoids analgous to the Weyl groups in the Bruhat decomposition for reductive algebraic groups. 

As the title might suggest, we will be considering not just Renner monoids, but a generalization of Renner monoids called Renner-Coxeter monoids. Renner-Coxeter monoids are to Renner monoids what Coxeter groups are to Weyl groups. As such, it will benefit us greatly to review some of the important results in the study of Coxeter groups. Many of these results can be found in \cite{Bjorner and Brenti} by Bjorner and Brenti, or are simple consequences of results from that source.

\begin{defn}
A {\bf Coxeter system} consists of a pair $(W, S)$ where

(i) $W$ is a group generated by $S$

(ii) there exists a function $m : S^{2} \rightarrow \{1, 2, \cdots, \infty\}$ so that,

\indent\indent\indent $m(s, t) = m(t, s)$ for all $s, t \in S$

\indent\indent\indent $m(s, t) = 1$ if and only if $s = t$

\indent\indent\indent $W$ has the presentation $\langle S \mid (st)^{m(s, t)}\text{ for all }s, t \in S \text{ with }m(s, t)\neq \infty\rangle$
\end{defn}

The function $m$ is often referred to as the {\bf Coxeter matrix}. The group, $W$, will be referred to as the {\bf Coxeter group} and $S$ will be termed the set of {\bf generators} or {\bf simple reflections}.

\begin{thm} \label{Coxeter length}
Define a function, $\ell : W \rightarrow \mathbb{N}$ called the {\bf length} function by, $\ell(w) = k$ where $k$ is minimum such that there exists an expression, $w = s_{1}s_{2} \cdots s_{k}$ with each $s_{i} \in S$. Any expression $w = s_{1}s_{2}\cdots s_{\ell(w)}$ with each $s_{i} \in S$ is called a {\bf reduced word} expression for $w$.

(1) $\ell(w^{-1}) = \ell(w)$

(2) $\ell(sw) = \ell(w) \pm 1$

(3) Two reduced words for $u$ must contain the exactly same generators, but possibly in different amounts and

\indent\indent arrangements.

(4) Any expression for an element $w \in W$, contains a reduced word expression for $w$ as a subword.

(5) $W$ is finite if and only if there is a unique element of maximum length, which is denoted $w_{0}$ and called the

\indent\indent{\bf longest element} of $W$.

(6) For any $u \in W$, $\ell(w_{0}u) = \ell(w_{0}) - \ell(u) = \ell(uw_{0})$

(7) $w_{0}^{2} = 1$
\end{thm}

\begin{proof}
(1) Proposition 1.4.2 (iv) in \cite{Bjorner and Brenti}.

(2) Proposition 1.4.2 (iii) in \cite{Bjorner and Brenti}.

(3) Corollary 1.4.8(ii) in \cite{Bjorner and Brenti}.

(4) Corollary 1.4.8(i) in \cite{Bjorner and Brenti}.

(5) Proposition 2.3.1 in \cite{Bjorner and Brenti}.

(6) Proposition 2.3.2 (ii) in \cite{Bjorner and Brenti}.

(7) Proposition 2.3.2 (i) in \cite{Bjorner and Brenti}.
\end{proof}

\begin{thm} \label{Coxeter Bruhat}
Let $u, v \in W$. We define the relation $u \leq v$ if and only if there is a reduced word expression, $v = s_{1}s_{2}\cdots s_{m}$ and a reduced word expression $u = s_{i_{1}}s_{i_{2}}\cdots s_{i_{n}}$ with $1 \leq i_{1} < i_{2} < \cdots < i_{n} \leq m$. This relation is known as the {\bf Bruhat order} and it has (among others) the follow properties.

(1) $\leq$ is a partial order on $W$.

(2) $u \leq v$ if and only if every reduced word for $v$ has a subword that is a reduced word for $u$.

(3) $u \leq v$ $\Rightarrow$ $\ell(u) \leq \ell(v)$

(4) $u \leq v$ $\Rightarrow$ $u^{-1} \leq v^{-1}$

(5) If $W$ is finite then $u \leq v$ $\Leftrightarrow$ $w_{0}v \leq w_{0}u$

(6) If $W$ is finite then $w \leq w_{0}$ for all $w \in W$.

(7) Suppose $u \leq v$ and $s \in S$. If $m = max\{v, sv\}$, $su \leq m$ 

(8) Suppose $u \leq v$ and $s \in S$. If $m = max\{v, vs\}$, $us \leq m$

(9) If $u < v$ there exists a chain $u = x_{0} < x_{1} < \cdots < x_{k} = v$ such that $\ell(x_{i}) = \ell(u) + i$ for all $i$.

(10) The {\bf Bruhat interval}, $[u, v] := \{w \in W \mid u \leq w \leq v\}$ is always finite. Indeed, $\abs{[u, v]} \leq 2^{\ell(v)}$.

(11) Suppose $u \leq v$, $x \leq y$ and $\ell(vy) = \ell(v) + \ell(y)$, then $ux \leq vy$
\end{thm}

\begin{proof}
(1) The properties of a partial order are clear from our definition.

(2) This is part of Corollary 2.2.3 in \cite{Bjorner and Brenti}.

(3) $\ell(u) = n \leq i_{n} \leq m = \ell(v)$

(4) Follows from the fact that if $s_{1}\cdots s_{k}$ is a reduced word for $w$ then $s_{k}\cdots s_{1}$ is a reduced word for $w$.

(5) (i) in Proposition 2.3.4 of \cite{Bjorner and Brenti}

(6) This is due to Proposition 2.3.1 in \cite{Bjorner and Brenti} and property (3) above.

(7) This is a quick consequence of the lifting property, Theorem 2.2.7 in \cite{Bjorner and Brenti}. The lifting property is a very important result which we will use in many proofs.

(8) Similar to (7).

(9) Theorem 2.2.6 in \cite{Bjorner and Brenti}.

(10) Corollary 2.2.4 in \cite{Bjorner and Brenti}.

(11) Since $u$ is a subword of $v$ and $x$ is a subword of $y$ it makes sense that a word for $ux$ is contained in a reduced word for $vy$. It follows that a reduced word for $ux$ is contained in a reduced word for $vy$.
\end{proof}

\begin{thm} \label{Coxeter cosets}
For a subset of generators, $I \subseteq S$ we denote by $W_{I}$ the subgroup of $W$ generated by $I$. A subgroup obtained in this way is called a {\bf standard parabolic subgroup}. There is a one to one correspondence between standard parabolic subgroups and subsets of $S$.

Let $I, J \subseteq S$, and let $w \in W$ be arbitrary. We define $^{I}W^{J} = \{u \in W \mid u \leq v, \forall v \in W_{I}uW_{J}\}$ as the set of minimum coset representatives for $I$ and $J$. Every double coset has a minimum element, and for $w$, we denote the minimum element of $W_{I}wW_{J}$ by $^{I}w^{J}$.

When $I = \emptyset$ we will often refer to $W^{J}$ and $w^{J}$, and when $J = \emptyset$ we will often refer to $^{I}W$ and $^{I}w$.

(1) $w \in $$^{I}W^{J}$ if and only if no reduced word for $w$ begins with a generator from $I$ or ends with a generator 

\indent\indent from $J$. 

(2) $^{I}w^{J} = $$^{I}(w^{J}) = ($$^{I}w)^{J}$

(3) For any $u, v \in W$, $u \leq v$ implies $^{I}u\phantom{~}^{J} \leq~$$^{I}v\phantom{~}^{J}$

(4) $^{I}W^{J} = $$^{I}W \cap W^{J}$

(5) For any $u \in W$, $\Big($$^{I}u\phantom{~\hspace{-1pt}}^{J}\Big)^{-1} =~$$^{I}\Big(u^{-1}\Big)\phantom{~}^{J}$

(6) For any $u, v \in W$, if $u = $$^{I}u^{J}$ then $u \leq v$ if and only if $u \leq w$ for all $w \in W_{I}vW_{J}$

(7) If $w \in $$^{I}W^{J}$ and $w^{\prime} \in W_{I}wW_{J}$ then there exist $u \in W_{I}$, $v \in W_{J}$ so that $w^{\prime} = uwv$ and $\ell(w^{\prime}) = \ell(u) + \ell(w) + \ell(v)$.

(8) Suppose that $u \leq v \in W$ and $v \in W^{I}$. If $w \in W_{I}$ then $uw \leq vw$.\label{multiplying by constant}

(9) Suppose that $u \leq v \in W$ and $v \in $$^{I}W$. If $w \in W_{I}$ then $wu \leq wv$.

(10) Suppose $a, b \in W^{I}$ and $x, y \in W_{I}$ with $ax \leq by$. Then we can find $u, v \in W_{I}$ so that $x = uv$, $\ell(x) = \ell(u) + \ell(v)$

\indent\indent and $au \leq b$, $v \leq y$.

(11) Suppose $a, b \in $$^{I}W$ and $x, y \in W_{I}$ with $xa \leq yb$. Then we can find $u, v \in W_{I}$ so that $x = vu$, $\ell(x) = \ell(v) + \ell(u)$

\indent\indent and $ua \leq b$, $v \leq y$.
\end{thm}

\begin{proof}
(1) Pick any $u \in W$. Pick any minimal element $v \in W_{I}uW_{J}$. Then it is clear that $v < sv$ and $v < vt$ for any $s \in I$, $t \in J$. So then no reduced word for $v$ can start with an element of $I$ or end with an element of $J$. Consider another element in $W_{I}uW_{J}$ with this property (it could be another minimal element, so this proof will also show us that each double coset does indeed have a minimum). Let this element be $v^{\prime} \in W_{I}uW_{J}$. 

We can find elements $x, x^{\prime} \in W_{I}$ and $y, y^{\prime} \in W_{J}$ so that $v = xuy$ and $v^{\prime} = x^{\prime}uy^{\prime}$. Let $x^{\prime}x^{-1} = a_{1}\cdots a_{m}$, $a_{i} \in I$ and $y^{-1}y^{\prime} = b_{1}\cdots b_{n}$, $b_{j} \in J$ be reduced words. Then $v^{\prime} = x^{\prime}uy^{\prime} = x^{\prime}x^{-1}vy^{-1}y^{\prime} = a_{1}\cdots a_{m}vb_{1}\cdots b_{n}$. By applying the lifting property $m + n$ many times we then see that $v^{\prime} \leq v$. By minimality of $v$, we have $v = v^{\prime}$. 

(2) Consider $x = w^{J}$ and $y = $$^{I}x$. By Lemma 2.4.3 in \cite{Bjorner and Brenti} we can find $u \in W_{I}$ so that $x = uy$ and $\ell(x) = \ell(u) = \ell(y)$. It is clear by (1) that no reduced word for $y$ starts with a generator from $I$. Suppose there is a reduced word for $y$ ending with $s \in J$. Then taking any reduced word for $u$ we can construct a reduced word for $x$ ending in $s$, a contradiction. So by (1) tells us that $^{I}w^{J} = y = $$^{I}x = $$^{I}(w^{J})$. Similarly, we can show that $^{I}w^{J} = ($$^{I}w)^{J}$. 

(3) It is clear that $^{I}u^{J} \leq u \leq v$. Since $^{I}v^{J} \in W_{I}vW_{J}$ we can find $a_{1}, \cdots, a_{m} \in I$ and $b_{1}, \cdots, b_{n} \in J$ so that $^{I}v^{J} = a_{m}\cdots a_{1}vb_{1}\cdots b_{n}$.

For each $i = 1, \cdots, m$, let $x_{i} = a_{i}x_{i-1}$, starting with $v = x_{0}$. Notice that $^{I}u^{J} \leq x_{0}$. Suppose $^{I}u^{J} \leq x_{i}$. Then either $x_{i+1} < x_{i}$ or $x_{i} < x_{i+1}$. In the former case, the lifting property guarantees that $^{I}u^{J} \leq x_{i+1}$ whereas in the latter, we see that $^{I}u^{J} \leq x_{i} < x_{i+1}$. So by induction $^{I}u^{J} \leq x_{m}$.

For each $i = 1, \cdots, m$, let $y_{j} = y_{j-1}b_{j}$, starting with $x_{m} = y_{0}$. Notice that $^{I}u^{J} \leq y_{0}$ by the preceeding paragraph. Suppose $^{I}u^{J} \leq y_{j}$. Then either $y_{j+1} < y_{j}$ or $y_{j} < y_{j+1}$. In the former case, the lifting property guarantees that $^{I}u^{J} \leq y_{j+1}$ whereas in the latter, we see that $^{I}u^{J} \leq y_{i} < y_{j+1}$. So by induction $^{I}u^{J} \leq y_{n} = $$^{I}v^{J}$.

(4) follows from (1) by considering reduced words.

(5) also follows from (1) by considering reduced words.

(6) $u \leq v$ implies $u = $$^{I}u^{J} \leq $$^{I}v^{J} \leq w$ for all $w \in W_{I}vW_{J}$ by definition of $^{I}v^{J}$. Conversely, if $u \leq w$ for all $w \in W_{I}vW_{J}$ then $u \leq v$.

(7) Lemma 2.4.3 in \cite{Bjorner and Brenti} notes that in the $W^{I}$ and $^{I}W$ cases this decomposition exists and is unique. Let $x \in W_{I}wW_{J}$. Then we can decompose $x = py$ with $\ell(x) = \ell(p) + \ell(y)$ and $y = $ $^{I}w$. Then let $y = qr$ with $\ell(y) = \ell(q) + \ell(r)$ and $q \in W^{J}$. By (2) then $q = $$^{I}x^{J}$ and $\ell(x) = \ell(p) + \ell(y) = \ell(p) + \ell(q) + \ell(r)$ as desired.

(8) Notice that $\ell(vw) = \ell(v) + \ell(w)$. As was remarked for the preceeding result, \cite{Bjorner and Brenti} notes that this decomposition is unique. So let $w = s_{1}\cdots s_{\ell(w)}$ be a reduced word decomposition for $w$. Then, since $u \leq v$ the resuilt follows by $\ell(w)$ applications of the lifting property.

(9) This is shown simiarly to (8).

(10) Write out $by = s_{1}\cdots s_{m}t_{1}\cdots t_{n}$ where $s_{1}\cdots s_{m} = b$ and $t_{1}\cdots t_{n} = y$. We know that since $ax \leq by$ we can find a reduced subword, $ax = s_{i_{1}} \cdots s_{i_{k}}t_{j_{1}}\cdots t_{j_{\ell}}$. Within this we can find a reduced word for $a$. But since $y \in W_{I}$ all $t_{i} \in I$ and since $a \in W^{I}$ no reduced word can end in $I$. Thus $a \leq s_{i_{1}} \cdots s_{i_{k}}$. Let $v = t_{j_{1}}\cdots t_{j_{\ell}}$ and let $u = a^{-1}axv^{-1}$. Then it is clear that $au = s_{i_{1}} \cdots s_{i_{k}} \leq b$ and $v \leq y$. Furthermore, by comparing lengths we can see that $x = uv$ with $\ell(x) = \ell(u) + \ell(v)$.

(11) Done similar to (10).
\end{proof}

\begin{thm} \label{Coxeter -fixes}
In Chapter 3 of \cite{Bjorner and Brenti} two new relations on $W$, related to the Bruhat order are defined. For $u, v \in W$ we say $u \leq_{L} v$ if and only if there exists $w \in W$ so that $v = wu$ with $\ell(v) = \ell(w) + \ell(u)$ and we say $u \leq_{R} v$ if and only if there exists $w \in W$ so that $v = uw$ with $\ell(v) = \ell(u) + \ell(w)$. These two relations are called the {\bf left} and {\bf right} {\bf weak orders} respectively.

(1) Both $\leq_{L}$ and $\leq_{R}$ are partial orders on $W$.

(2) $u \leq_{R} v$ $\Leftrightarrow$ $\ell(v) = \ell(u) + \ell(u^{-1}v)$ and $u \leq_{L} v$ $\Leftrightarrow$ $\ell(v) = \ell(u) + \ell(vu^{-1})$

(3) $u \leq_{L} v$ $\Rightarrow$ $u \leq v$, $u \leq_{R} v$ $\Rightarrow$ $u \leq v$

(4) Suppose that $W$ is finite. Then for any $w \in W$, $w \leq_{L} w_{0}$ and $w \leq_{R} w_{0}$

(5) Suppose $u \leq v$, $y \leq x$ and $x^{-1} \leq_{L} u$. Then $ux \leq vy$

(6) Suppose $v \leq u$, $x \leq y$ and $u^{-1} \leq_{R} x$. Then $ux \leq vy$

(7) Take $I \subseteq S$ and suppose $W$ is finite. For any $u \in W^{I}$, $v \in  W_{w_{0}Iw_{0}}w_{0}W_{I}$, $u \leq_{L} v$

(8) Take $I \subseteq S$ and suppose $W$ is finite. For any $u \in $$^{I}W$, $v \in  W_{I}w_{0}W_{w_{0}Iw_{0}}$, $u \leq_{R} v$
\end{thm}

\begin{proof}
(1) Clear by the definition of our weak orders.

(2) (ii) in Proposition 3.1.2 of \cite{Bjorner and Brenti}.

(3) Is clear by definition of $\leq$ as $u$ is definitely a subword of $v$.

(4) (iii) in Proposition 3.1.2 of \cite{Bjorner and Brenti}.

(5) Consider $x = s_{1}s_{2}\cdots s_{\ell}$ as a reduced word. We will prove the result by induction on $\ell$. For $\ell = 0$, $x = y = 1$ and we already know that $u \leq v$, so $ux \leq vy$. Assume the result is true for all $\ell \leq k$ and suppose that $\ell(x) = k+1$. Notice that $s_{1}s_{2}\cdots s_{k} = xs_{k+1} < x$ and also $(xs_{k+1})^{-1} \leq_{L} u$ as well. Consider $ys_{k+1}$. By the lifting property we have either $y \leq xs_{k+1}$ or $y \not\leq xs_{k+1}$ and $ys_{k+1} \leq xs_{k+1}$.

If $y \leq xs_{k+1}$ then $u(xs_{k+1}) \leq vy$ by induction and it follows that $ux = u(xs_{k+1})s_{k+1} < u(xs_{k+1}) \leq vy$. If $y \not\leq xs_{k+1}$ and $ys_{k+1} \leq xs_{k+1}$ then by induction $u(xs_{k+1}) \leq v(ys_{k+1})$. By applying the lifting property it can be quickly shown that $ux = u(xs_{k+1})s_{k+1} \leq v(ys_{k+1})s_{k+1} = vy$.

(6) Similar to (5).

(7) First notice that for any $s \in S$ $\ell(w_{0}sw_{0}) = \ell(w_{0}) - \ell(sw_{0}) = \ell(w_{0}) - (\ell(w_{0}) - \ell(s)) = \ell(s) = 1$, so $w_{0}sw_{0}$ is a generator. Thus, $w_{0}Iw_{0} \subseteq S$ and so our statement makes sense. 

Consider any element $v = aw_{0}b \in W_{w_{0}Iw_{0}}w_{0}W_{I}$. Then we can write $a = s_{1}\cdots s_{k}$ and $b = t_{1}\cdots t_{\ell}$ as reduced words. So, $aw_{0}b = s_{1}s_{2}\cdots s_{k}w_{0}t_{1}t_{2}\cdots t_{\ell} = w_{0}s^{\prime}_{1}w_{0}w_{0}s^{\prime}_{2}w_{0}\cdots w_{0}s^{\prime}_{k}w_{0}w_{0}t_{1}t_{2}\cdots t_{\ell} = w_{0}s^{\prime}_{1}s^{\prime}_{2}\cdots s^{\prime}_{k}t_{1}t_{2}\cdots t_{\ell}$ with all the $s^{\prime}_{i}$ and $t_{k}$ being simple reflections in $I$. So $aw_{0}b \in w_{0}W_{I}$ and indeed, we can even show $W_{w_{0}Iw_{0}}w_{0} = W_{w_{0}Iw_{0}}w_{0}W_{I} = w_{0}W_{I}$ by symmetry. It follows that $^{w_{0}Iw_{0}}w_{0} = $$^{w_{0}Iw_{0}}w_{0}^{I} = w_{0}^{I}$.

Let $u \in W^{I}$. By (4) we can see that $uw_{0}(I) \leq_{L} w_{0} = w_{0}^{I}w_{0}(I)$. By Proposition 3.1.2(vi) in \cite{Bjorner and Brenti} it follows that $u \leq_{L} w_{0}^{I}$. For any $v \in W_{w_{0}Iw_{0}}w_{0}W_{I}$ we can write $v = x$$^{w_{0}Iw_{0}}w_{0}$, so $w_{0}^{I} = $$^{w_{0}Iw_{0}}w_{0} \leq_{L} v$. Since $\leq_{L}$ is a partial order, transitivity lets us conclude that $u \leq_{L} v$.

(8) Similar to (7).
\end{proof}

These results are meant to highlight the important properties from Coxeter group theory which will permeate all our later results. There are other results we will need, which will be invoked at the time. Those wishing a much more complete understanding of the Bruhat order and minimum coset representatives are encouraged to read Chapters 1 and 2 of \cite{Bjorner and Brenti}.

\section{Renner-Coxeter Systems}

In his  paper, \cite{Godelle}, Eddy Godelle made an excellent definition that went beyond the usual Renner monoid derived from reductive monoids in algebraic geometry. Taking the basic structures of the Renner monoid from reductive monoids he created a whole collection of factorizable monoids with similar, simple, and robust structure.

A factorizable monoid is one which is unit regular ($M = E(M)G(M) = G(M)E(M)$) and its set of idempotents is a meet semilattice with the partial order $e \leq f \Leftrightarrow ef = fe = e$. Factorizable monoids are known to be inverse monoids. That is, every element $m \in M$ has a unique element, $n$, called the inverse of $m$, so that $mnm = m$ and $nmn = n$.

\begin{defn}
A {\bf generalized Renner-Coxeter system} is a triple $(\Ren, \Lambda, S)$ such that,

(i) $\Ren$ is a factorizable monoid, whose group of units we will denote by $W$

(ii) $\Lambda$ is a transversal of $E(\Ren)$ for the action of $W$ and a sub-meet semilattice

(iii) $(W, S)$ is a Coxeter system

(iv) for every pair $e \leq e^{\prime} \in E(\Ren)$, there exists $w \in W$ and $f \leq f^{\prime} \in \Lambda$ so that $wew^{-1} = f$

\hspace{20pt} and $we^{\prime}w^{-1} = f^{\prime}$

(v) the map $\lambda^{*} : \Lambda \rightarrow P(S)$ given by $\lambda^{*}(e) = \{s \in S \mid se = es \neq e\}$ satisfies $e\leq f$ $\Rightarrow$ $\lambda^{*}(e) \subseteq \lambda^{*}(f)$

(vi) for every $e \in \Lambda$, the groups, $\{w \in W \mid we = ew\}$ and $\{w \in W \mid we = ew = e\}$, are standard

\hspace{20pt} parabolic subgroups of $W$
\end{defn}

In Godelle's paper he refers to the monoid, $\Ren$ as a {\emph generalized Renner monoid}, however in this paper we will refer to $\Ren$ as a {\bf Renner-Coxeter monoid} as it is a little more clear at a glance that we are talking about a monoid described by the above system rather than the more familiar Renner monoids from reductive monoids with which we will occasionally contrast. We shall call $\Lambda$ the {\bf cross-sectional lattice} and $S$ again as the set of {\bf generators} or the set of {\bf simple reflections}.

For the remainder of this paper, $\Ren$ will always represent a Renner-Coxeter monoid, $W$ will always be its group of units, $S$ its generators, $\Lambda$ the cross-sectional lattice, and so on. Any other assumptions on these basic structures will be stated during results that require them.

In \cite{Godelle}, Godelle covers three different examples showing that Renner monoids derived from reductive algebraic monoids (Example 1.7 in \cite{Godelle}), abstract monoids of finite Lie type (Example 1.8 in \cite{Godelle}), and Kac-Moody groups over fields of characteristic $0$ (Example 1.9 in \cite{Godelle}) are all Renner-Coxeter systems. This allows our work to cover much more than just reductive monoid results, although later we will make an assumption of finiteness aimed at examining those Renner monoids more specifically.

As they are parabolic subgroups of $W$, we can define maps $\lambda : \Lambda \rightarrow P(S)$ and $\lambda_{*} : \Lambda \rightarrow P(S)$ so that $\{w \in W \mid we = ew\} = W_{\lambda(e)}$ and $\{w \in W \mid we = ew = e\} = W_{\lambda_{*}(e)}$ for each $e \in \Lambda$.  We term $\lambda$, $\lambda_{*}$, and $\lambda^{*}$ as the {\bf type map}, {\bf lower type map}, and {\bf upper type map} respectively. These type maps have a special relationship which will be explored by the next two results.

\begin{lem} \label{commuting parabolic}
For $J \subseteq I \subseteq S$ the following are equivalent

(1) $J$ and $I\remove J$ commute

(2) for any $u \in W_{I\remove J}$ and $v \in W_{J}$, $uv = vu$

(3) $W_{J}$ is a normal subgroup of $W_{I}$

(4) $W_{I\remove J}$ is a normal subgroup of $W_{I}$

(5) $W_{I} = W_{I\remove J}\times W_{J}$

(6) $W_{I} = W_{J}\times W_{I\remove J}$
\end{lem}

\begin{proof}
We shall just show (1) $\Rightarrow$ (2) $\Rightarrow$ (3) $\Rightarrow$ (5) $\Rightarrow$ (1) and note that (1) $\Rightarrow$ (2) $\Rightarrow$ (4) $\Rightarrow$ (6) $\Rightarrow$ (1) is done similarly, completing the result.

(1) $\Rightarrow$ (2), (2) $\Rightarrow$ (3) are clear. For (3) $\Rightarrow$ (5) let $w \in W_{I}$ with a given reduced word $w = s_{1}s_{2}\cdots s_{\ell}$. Define $a_{0} = b_{0} = 1$. We will build up two elements $a_{\ell} \in W_{I\remove J}$ and $b_{\ell} \in W_{J}$ so that $w = a_{\ell}b_{\ell}$. At step $i \geq 1$, if $s_{i} \in J$ then let $a_{i} = a_{i-1} \in W_{I\remove J}$ and let $b_{i} = b_{i-1}s_{i} \in W_{J}$. If instead, $s_{i} \in I\remove J$ then let $a_{i} = a_{i-1}s_{i}$ and let $b_{i} = s_{i}b_{i-1}s_{i} \in W_{J}$ by our normality assumption.

Observe that at each stage, $a_{i} \in W_{I\remove J}$, $b_{i} \in W_{J}$ and $w = a_{i}b_{i}s_{i+1}s_{i+2}\cdots s_{\ell}$. Thus, after $\ell$ steps of the process, we have $w = a_{\ell}b_{\ell}$. Since $w$ was arbitrary we see that $W_{I} = W_{I\remove J}\times W_{J}$.

For (5) $\Rightarrow$ (1) take $i \in I\remove J$ and $j \in J$. Either $iji = i$, $iji = j$, or $iji$ is a reduced word. The first case implies $i = j$, a contradiction. The second can rearrange to show $ij = ji$. Our result will be proven if we can just show that $iji$ is not a reduced word.

Clearly, $iji \in W_{I}$, so $iji = uv$ for some $u \in W_{I\remove J}$ and $v \in W_{J}$. Let $w$ be a reduced word subword (recall Theorem \ref{Coxeter length}) for $uv$. Then $w = s_{1}\cdots s_{k}s_{k+1}\cdots s_{\ell(w)}$ where $k$ is such that $s_{a} \in I\remove J$ for all $a \leq k$ and $s_{a} \in J$ for all $k < a$. But two reduced words for the same group element must have the exact same $s \in S$, just in different amounts (Theorem \ref{Coxeter length}) so $s_{a} = i$ for $a \leq k$, $s_{a} = j$ for $a > k$. Thus, $iji = w = i^{k}j^{\ell(w) - k}$. Since $\ell(w) = 3$ we either have $iji = iij = j$ or $iji = ijj = i$, both of which show that $iji$ is not a  reduced word.
\end{proof}

\begin{prop} \label{lower type map is normal}
For any $e \in \Lambda$, $W_{\lambda_{*}(e)}$ is a normal subgroup of $W_{\lambda(e)}$.
\end{prop}

\begin{proof}
Let $u \in W_{\lambda_{*}(e)}$ and $v \in W_{\lambda(e)}$. By definition we know that $ue = e = eu$. Consider $vuv^{-1}$ and notice that $vuv^{-1}e = vuev^{-1} = vev^{-1} = evv^{-1} = e = vv^{-1}e = vev^{-1} = veuv^{-1} = evuv^{-1}$. Thus $vuv^{-1} \in W_{\lambda_{*}(e)}$ as desired.
\end{proof}

Since $W_{\lambda_{*}(e)} \unlhd W_{\lambda(e)}$ we can now use any of the six properties of the preceeding lemma to our advantage. We will make great use later on of facts like $\lambda_{*}(e)$ and $\lambda^{*}(e)$ commute and $W_{\lambda(e)} = W_{\lambda_{*}(e)}\times W_{\lambda^{*}(e)} = W_{\lambda^{*}(e)}\times W_{\lambda_{*}(e)}$.

In addition to introducing the concept of generalized Renner-Coxeter systems, \cite{Godelle} lets us understand the presentation of such systems, introduces a notion of length, and also gives us the following result, of which we will make great use. 

\begin{thm}
Suppose that $e, f \in \Lambda$ and $w \in $$^{\lambda(e)}W^{\lambda(f)}$. There is a unique element, $g \in \Lambda$ such that,

(1) $w \in W_{\lambda_{*}(g)}$

(2) $ewf = g = fw^{-1}e$

(3) $g = max\{h \in \Lambda \mid h \leq e, h\leq f, w \in W_{\lambda(h)}\}$

\noindent This element is denoted, $e\wedge_{w}f$, and shall be referred to as {\bf Godelle's meet}.
\end{thm}

\begin{proof}
This is actually two results from \cite{Godelle} written together. (1) and $ewf = g$ follow from Lemma 1.12 and (3) and $g = fw^{-1}e$ are stated as part (iii) of Corollary 1.13.
\end{proof}

Our last piece of introductory monoid material is the following collection of equivalence relations, familiar to nearly anyone who has studied semigroups before.

\begin{defn}
For any $r, s \in \Ren$, we recall the following defintions for {\bf Green's relations}.

(1) $r \Jc s$ if and only if $\Ren r\Ren = \Ren s\Ren$ if and only if there are $u, v \in W$ so that $s = urv$

(2) $r \Lc s$ if and only if $\Ren r = \Ren s$ if and only if there is $u \in W$ so that $s = ur$

(3) $r \Rc s$ if and only if $r\Ren = s\Ren$ if and only if there is $v \in W$ so that $s = rv$

(4) $r \Hc s$ if and only if $r \Lc s$ and $r \Rc s$ if and only if there are $e, f \in E(\Ren)$ and $u, v \in W$ so that $r = eu = uf$ 

\indent\indent and $s = ev = vf$
\end{defn}

That these are equivalence relations on a monoid is a well known result. Of particular interest are the definitions for the relations in terms of elements of $W$ which come about since $\Ren$ is by definition a factorizable monoid and hence a unit regular monoid.

Semigroup theorists will notice that the common $\Dc$ relation is absent in our definition. However, upon closer inspection one can see that $\Jc = \Lc\circ\Rc$ which is the definition of the $\Dc$ relation.

As we move forward we will be seeing a lot of Green's relations. When we wish to denote a generic relation, we shall often use the $\Tc$ symbol. To denote a specfic equivalence class, we shall use $T_{r} = \{s \in \Ren \mid r\Tc s\}$ to mean the $\Tc$-class of the element $r$.

\section{Adherence Orders on Renner-Coxeter Monoids}

Our goal is to generalize the Adherence order of Renner monoids from reductive algebraic monoids (\cite{Renner Bruhat I}). Thanks to excellent work done by Pennell, Putcha, and Renner in \cite{PPR} the Adherence order between two elements of the Renner monoid can be determined by first decomposing the elements in question in a unique way, known as the standard form decomposition. It turns out (as has been noted by Godelle in his paper) the elements of a Renner-Coxeter monoid can be decomposed similarly. We will use these standard forms in the definition of our Adherence orders.

\begin{prop}
Let $r \in \Ren$.

(1) There exists a unique set of elements, $x, e, y \in \Ren$ so that $r = xey$, $e \in \Lambda$, $x \in W^{\lambda_{*}(e)}$, and $y \in $$^{\lambda(e)}W$. When we 

\indent\indent write $r = xey$ in this way we say $r$ is written in {\bf left standard form}.

(2) There exists a unique set of elements, $y, e, x \in \Ren$ so that $r = yex$, $e \in \Lambda$, $x \in $$^{\lambda_{*}(e)}W$, and $y \in W^{\lambda(e)}$. When we 

\indent\indent write $r = yex$ in this way we say $r$ is written in {\bf right standard form}.

(3) There exists a unique set of elements, $x, e, y, z \in \Ren$ so that $r = xeyez$, $e \in \Lambda$, $x \in W^{\lambda(e)}$, $y \in W_{\lambda^{*}(e)}$, and

\indent\indent  $z \in $$^{\lambda(e)}W$. When we write $r = xeyez$ in this way we say $r$ is written in {\bf hybrid standard form}.

(4) If $r = xeyez$ is in hybrid standard form, then $r = (xy)ez$ is in left standard form and $r = xe(yz)$ is in right

\indent\indent  standard form.
\end{prop}

\begin{proof}
(1) and (2) are Proposition 1.11 in \cite{Godelle}.

(3) and (4) Let us take $r = aeb$ to be in left standard form. Since $a \in W^{\lambda_{*}(e)}$ and $W_{\lambda(e)} = W_{\lambda^{*}(e)}\times W_{\lambda_{*}(e)}$ is is not hard to see that if $x = a^{\lambda(e)}$ then there exists $y \in W_{\lambda^{*}(e)}$ so that $a = xy$. Letting $z = b \in $$^{\lambda(e)}W$ we can see that $r = aeb = xyez = xeyez$ satisfies the properties of hybrid standard form.

Suppose $r = xeyez$ satisfies the conditions of hybrid standard form. Since $e$ is an idempotent and $y \in W_{\lambda^{*}(e)}$ we see that $ye = yee = eye = eey = ey$. So it is clear that $r = xeyez = (xy)ez = xe(yz)$. Since $x \in W^{\lambda(e)}$ and $y \in W_{\lambda^{*}(e)}$ we can see that $xy \in W^{\lambda_{*}(e)}$, so $(xy)ez$ is in left standard form and similarly $xe(yz)$ is in right standard form.

To show uniqueness, suppose that $r = xeyez = aebec$ are two such decompositions. Then $(xy)ez = (ab)ec$ are in left standard form, which is unique. Thus $z = c$ and $xy = ab$. Then $x = (xy)^{\lambda(e)} = (ab)^{\lambda(e)} = a$ and it follows quickly that $y = b$ as well.
\end{proof}

With our standard forms in hand, we can more easily compare some of Green's relations by observation of the decomposition.

\begin{prop} \label{standard form Green's}
Suppose $r = xey$ and $s = afb$ are elements of $\Ren$ written in left standard form. Then, 

(1) $r\Jc s$ if and only if $e = f$

(2) $r\Lc s$ if and only if $e = f$ and $y = b$

\noindent Suppose $r = xey$ and $s = afb$ are elements of $\Ren$ written in right standard form. Then, 

(3) $r\Jc s$ if and only if $e = f$

(4) $r\Rc s$ if and only if $e = f$ and $y = b$

\noindent Suppose $r = aebec$ and $s = xfyfz$ are elements of $\Ren$ written in hybrid standard form. Then, 

(5) $r\Jc s$ if and only if $e = f$

(6) $r\Lc s$ if and only if $e = f$ and $c = z$

(7) $r\Rc s$ if and only if $e = f$ and $a = x$

(8) $r\Hc s$ if and only if $e = f$, $a = x$ and $c = z$
\end{prop}

\begin{proof}
(1) $r \Jc s$ if and only if $WrW = WsW$ by definition. Then we see that $WeW = WrW = WsW = WfW$, so $r \Jc s$ if and only if $e = f$.

(2) $r \Lc s$ if and only if $Wr = Ws$ by definition. So $r \Lc s$ if and only if we can find $w \in W$ so that $r = ws$ or rather, $xey = wafb = (wa)^{\lambda_{*}(f)}fb$. The first and last of these are in left standard form which is unique. So $x = (wa)^{\lambda_{*}(f)}$, $e = f$, and $y = b$. The converse follows quickly.

(3) and (4) are done similarly to (1) and (2).

(5), (6), and (7) are each covered by previous results thanks to (4) from the preceeding proposition. 

(8) follows from applying (6) and (7) and the fact that $r \Hc s$ if and only if $r \Lc s$ and $r \Rc s$.
\end{proof}

(1) in the preceeding shows us that $\Lambda$ really is cross-sectional with respect to the $\Jc$-classes of $\Ren$, earning its name.

\begin{defn}
Let us define two distinct relations on $\Ren$ based on our two standard forms. Let $r, s \in \Ren$ be arbitrary elements.

(1) If $r = xey$ and $s = afb$ are in left standard form, we say $r \leq^{+} s$ if and only if $e \leq f$ and there exists an 

\indent\indent element $w \in W_{\lambda^{*}(f)}W_{\lambda_{*}(e)}$ so that $x \leq aw$ and $w^{-1}b \leq y$.

(2) If $r = yex$ and $s = bfa$ are in right standard form, we say $r \leq^{-} s$ if and only if $e \leq f$ and there exists an

\indent\indent  element $w \in W_{\lambda_{*}(e)}W_{\lambda^{*}(f)}$ so that $x \leq wa$ and $bw^{-1} \leq y$.

\noindent These partial orders are referred to as the {\bf Adherence orders} on $(\Ren, \Lambda, S)$.
\end{defn}

Our choice of associating $+$ to the left standard form comes from the historic focus on left standard form (\cite{PPR}). We will show shortly that the notation, $\leq^{+}$ and $\leq^{-}$ is apt, as these are partial orders. But towards the proof we first prove the following lemma.

\begin{lem} \label{circ max exits}
For any $u, v \in W$, the following exist and are equal

(1) $max\{u^{\prime}v \mid u^{\prime} \leq u\}$

(2) $max\{uv^{\prime} \mid v^{\prime} \leq v\}$

(3) $max\{u^{\prime}v^{\prime} \mid u^{\prime} \leq u, v^{\prime} \leq v\}$

\noindent As a result, we denote this element by $u\circ v$ and, as it produces a maximum with respect to the Bruhat order, it is called the {\bf optimization operator}.
\end{lem}

\begin{proof}
We are going to show that (1) and (3) exist and are equal by induction on $\ell(u)$. The existence and equality of (2) and (3) follow by a similar argument, completing the result.

Let us induct on $k = \ell(u)$. If $k = 0$ then $u = 1$ and the result is clear. Suppose that $k > 0$. Then we can find $s \in S$ so that $w := su < u$. By induction, $max\{w^{\prime}v \mid w^{\prime} \leq w\}$ and   $max\{w^{\prime}v^{\prime} \mid w^{\prime} \leq w, v^{\prime} \leq v\}$ exist and are equal. Let us call them $x$. Let $m = max\{x, sx\}$. Now it suffices to show that $m$ is the maximum of each of the sets specified in (1) and (3).

(1) We can find $w_{1} \leq w$ so that $x = w_{1}v$, so $m = w_{1}v$ or $sw_{1}v$ and since $w_{1}, sw_{1} \leq max\{w, sw\} = u$ we see that $m \in \{u^{\prime}v \mid u^{\prime} \leq u\}$. Take any $u^{\prime} \leq u$. If $su^{\prime} < u^{\prime}$, then $su^{\prime} \leq w$, in which case $su^{\prime}v \leq x$ and hence $u^{\prime}v = ssu^{\prime}v \leq m$. On the other hand, if $u^{\prime} < su^{\prime}$, then $u^{\prime} \leq w$ (by the lifting property), in which case $u^{\prime}v \leq x$ and hence $u^{\prime}v \leq m$.

(3) We can find $w_{1} \leq w$ and $v_{1} \leq v$ so that $x = w_{1}v_{1}$, so either $m = w_{1}v_{1}$ or $sw_{1}v_{1}$ and since $w_{1}, sw_{1} \leq max\{w, sw\} \leq u$ we see that $m \in \{u^{\prime}v^{\prime} \mid u^{\prime} \leq u, v^{\prime} \leq v\}$. Take any $u^{\prime} \leq u$, $v^{\prime} \leq v$. If $su^{\prime} < u^{\prime}$, then $su^{\prime} \leq w$, in which case $su^{\prime}v^{\prime} \leq x$ and hence $u^{\prime}v^{\prime} = ssu^{\prime}v^{\prime} \leq m$. On the other hand, if $u^{\prime} < su^{\prime}$, then $u^{\prime} \leq w$ (by the lifting property), in which case $u^{\prime}v^{\prime} \leq x$ and hence $u^{\prime}v^{\prime} \leq m$.
\end{proof}

We denote the element of $W$ in the preceeding lemma by $u\circ v$ just as Putcha does in \cite{Putcha Shellability}. $\circ$ turns out to be a very interesting and important binary operation. Of particular note, $uv = u\circ v$ if and only if $\ell(uv) = \ell(u) + \ell(v)$ for every $u, v \in W$.

Not only are both $\leq^{+}$ and $\leq^{-}$ partial orders on $\Ren$, but they can be seen to quite naturally extend some of the other, more familiar partial orders on important subsets of $\Ren$. Both $\leq^{+}$ and $\leq^{-}$ extend the Bruhat order on the group of units, $W$, as well as extend the usual product partial order for idempotents, $E(\Ren)$.

\begin{thm} \label{The Partial Order Theorem}
For $\varepsilon = +\text{ or }-$,

(1) $(\Ren, \leq^{\varepsilon})$ is a partially ordered set

(2) $(W, \leq)$ is a sub-partially ordered set of $(\Ren, \leq^{\varepsilon})$

(3) $(\Lambda, \leq)$ is a sub-partially ordered set of $(\Ren, \leq^{\varepsilon})$

(4) $(E(\Ren), \leq)$ is a sub-partially ordered set of $(\Ren, \leq^{\varepsilon})$
\end{thm}

\begin{proof}
(1) The result for $\leq^{+}$ and $\leq^{-}$ are symmetrical, so we will just demonstate it for $\varepsilon = +$. We start with reflexivity. Since $e \leq e$ and $1 \in  W_{\lambda^{*}(e)}W_{\lambda_{*}(e)}$ it is clear that $xey \leq^{+} xey$. Suppose $xey \leq afb$ and $afb \leq xey$. Then $e \leq f \leq e$ which implies $e = f$. We can find $u, v \in W_{\lambda(e)}$ so that $x\leq au$, $a \leq xv$, $u^{-1}b \leq y$, and $v^{-1}y \leq b$. But since $y, b \in $$^{\lambda(e)}W$, $b \leq u^{-1}b \leq y$ and $y \leq v^{-1}y \leq b$. Thus $y = b$, and $u = v = 1$. From there we can see that $x \leq a \leq x$ or rather $x = a$. Hence $xey = afb$.

Now, towards transitivity, suppose that $xey \leq^{+} afb$ and $afb \leq^{+} cgd$. It is clear that $e \leq g$. We can find $u \in W_{\lambda^{*}(f)}$, $p \in W_{\lambda_{*}(e)}$ and $v \in W_{\lambda^{*}(g)}$, $q \in W_{\lambda_{*}(f)}$ so that $x\leq aup$, $a \leq cvq$, $p^{-1}u^{-1}b \leq y$, and $q^{-1}v^{-1}d \leq b$. By Theorem \ref{Coxeter cosets}, we see that $x\leq aup$, $a \leq cvq$ imply $x\leq au$, $a \leq cv$.

Since $x \leq au$ the subword property tells us that $x = a_{1}u_{1}$ for some $a_{1} \leq a$ and $u_{1} \leq u$. Then we can see that $a_{1} \leq a \leq cv$ and hence $x = a_{1}u_{1} \leq (cv)\circ u_{1} = (cv)u_{2}$ for some $u_{2} \leq u_{1} \leq u \in W_{\lambda^{*}(f)} \subseteq W_{\lambda^{*}(g)}$ (by Lemma \ref{circ max exits}). Then by Theorem \ref{multiplying by constant}, $q^{-1}u_{2}^{-1}v^{-1}d = u_{2}^{-1}q^{-1}v^{-1}d \leq u_{2}^{-1}b \leq u^{-1}b$. We know that $p^{-1}u^{-1}b \leq y$ so it follows that $^{\lambda_{*}(e)}(q^{-1}u_{2}^{-1}v^{-1}d) \leq p^{-1}u^{-1}b \leq y$. This means that we can find an element $r \in W_{\lambda_{*}(e)}$ such that $^{\lambda_{*}(e)}(q^{-1}u_{2}^{-1}v^{-1}d) = r^{-1}q^{-1}u_{2}^{-1}v^{-1}d$. 

Notice that $qr \in W_{\lambda_{*}(f)}W_{\lambda_{*}(e)} \subseteq W_{\lambda_{*}(e)}$ and $vu_{2} \in W_{\lambda^{*}(g)}$. By seeing $vu_{2}qr \in W_{\lambda^{*}(g)}W_{\lambda_{*}(e)}$ and $x \leq cvu_{2}qr$ and $r^{-1}q^{-1}u_{2}^{-1}v^{-1}d \leq y$ we may conclude $xey \leq^{+} cgd$.

(2) For $\varepsilon = +$, notice that for any $u, v \in W$, $u\cdot 1\cdot 1$ and $v\cdot 1\cdot 1$ are the left standard forms. So we see that $u\leq^{+} v$ if and only if we can find $w \in W_{\lambda^{*}(1)}W_{\lambda_{*}(1)} = W$ so that $u \leq vw$ and $w^{-1}1 \leq 1$. The latter condition forces $w = 1$ and so it is equivalent to $u \leq v$ in the usual Bruhat order. A similar proof can be given for $\leq^{-}$.

(3) Take $e, f \in \Lambda$. Observe that $1\cdot e\cdot 1$, $1\cdot f\cdot 1$ are the respective standard forms (regardless of $\varepsilon$). Then, depending on choice of $\varepsilon$, we see that $e \leq^{\varepsilon} f$ if and only if $e \leq f$ and there exists $w \in W_{\lambda^{*}(f)}W_{\lambda_{*}(e)}$ such that $1\leq w$ and $w^{-1} \leq 1$ or $e \leq f$ and there exists $w \in W_{\lambda_{*}(e)}W_{\lambda^{*}(f)}$ such that $1\leq w$ and $w^{-1} \leq 1$. In either case, $w = 1$ and we are left with $e \leq f$.

(4) We shall just prove the result for $\varepsilon = +$, as $\varepsilon = -$ is done in a similar manner. Let $ueu^{-1}, vfv^{-1}$ be arbitrary idempotents, with $e, f \in \Lambda$. One can quickly check that it is safe to assume $u \in W^{\lambda(e)}$ and $v \in W^{\lambda(f)}$. Suppose that $ueu^{-1} \leq^{+} vfv^{-1}$. Then $e \leq f$ and we can find $w \in W_{\lambda^{*}(f)}$ and $x \in W_{\lambda_{*}(e)}$ so that $vwx \leq u \leq vw$. Thus $u \in vwW_{\lambda_{*}(e)}$, and we can find $y \in W_{\lambda_{*}(e)}$ so that $u = vwy$. We can rearrange and see that $u^{-1}v = y^{-1}w^{-1}$ and $v^{-1}u = wy$. Observe that $(ueu^{-1})(vfv^{-1}) = uey^{-1}w^{-1}fv^{-1} = uefw^{-1}v^{-1} = ueu^{-1}$ and $(vfv^{-1})(ueu^{-1}) = vfwyeu^{-1} = vwfeu^{-1} = ueu^{-1}$.

Suppose $(ueu^{-1})(vfv^{-1}) = ueu^{-1} = (vfv^{-1})(ueu^{-1})$. We write $v^{-1}u = pqr$ where $p \in W_{\lambda(f)}$, $q \in $$^{\lambda(f)}W^{\lambda(e)}$, and $r \in W_{\lambda(e)}$. Thus $ueu^{-1} = (vfv^{-1})(ueu^{-1}) = vpfqeru^{-1} = vp(f\wedge_{q}e)ru^{-1} \in W(f\wedge_{q}e)W$. Since $f\wedge_{q}e \in \Lambda$ it follows that $f\wedge_{q}e = e$, $q \in W_{\lambda_{*}(e)}$ and thus $f\wedge_{q}e = fqe = fe = e$. In a similar fashion we can show that $ef = e$. So we have demonstrated that $e \leq f$. In fact, we also know $ueu^{-1} = vfv^{-1}ueu^{-1} = vpfqeru^{-1} = vperu^{-1} = vpreu^{-1}$. It follows that $vpr \in uW_{\lambda_{*}(e)}$. This means we can find $z \in W_{\lambda_{*}(e)}$ so that $vprz = u$. By using the commuting properties of $W_{\lambda_{*}(f)} \subseteq W_{\lambda_{*}(e)}$ and $W_{\lambda^{*}(e)} \subseteq W_{\lambda^{*}(f)}$ we can rearrange this equation and find terms $a \in W_{\lambda^{*}(f)}$ and $b \in W_{\lambda_{*}(e)}$ so that $u = vprz = vab$. From there it follows that $u \leq v(ab)$ and $(ab)^{-1}v^{-1} \leq u^{-1}$ as desired.
\end{proof}

We are going to primarily be concerning ourselves with $\leq^{+}$, but nearly everything we discover will apply (with some symmetry) to $\leq^{-}$ and the right standard form. When possible, we will make joint statements using the above $\varepsilon$ notation, however, if that would become too cumbersome, we will state the result in terms of one of the two (mostly $\leq^{+}$) and leave curious readers to mirror the arguments and come up with the analogous statement.

\begin{prop} \label{order and *}
For any two elements, $r, s\in \Ren$, the following are equivalent.

(1) $r \leq^{+} s$ in $(\Ren, \Lambda, S)$

(2) $r^{*} \leq^{-} s^{*}$ in $(\Ren, \Lambda, S)$
\end{prop}

\begin{proof}
Let $r = xey$ and $s = afb$ be written in left standard form. It is quick to see that $r^{*} = y^{-1}ex^{-1}$ and $s^{*} = b^{-1}fa^{-1}$ are written in right standard form. Then we note the following equivalences, 

\vspace{5pt}\begin{tabular}{rcl}
$r \leq^{+} s$ & if and only if & $e \leq f$ and there exists $w \in W_{\lambda^{*}(f)}W_{\lambda_{*}(e)}$ so that $x \leq aw$ and $w^{-1}b \leq y$ \\

~ & if and only if & $e \leq f$ and there exists $w \in W_{\lambda^{*}(f)}W_{\lambda_{*}(e)}$ so that $x^{-1} \leq w^{-1}a^{-1}$ and $b^{-1}w \leq y^{-1}$ \\

~ & if and only if & $e \leq f$ and there exists $w \in W_{\lambda_{*}(e)}W_{\lambda^{*}(f)}$ so that $x^{-1} \leq wa^{-1}$ and $b^{-1}w^{-1} \leq y^{-1}$ \\

~ & if and only if & $r^{*} \leq^{-} s^{*}$
\end{tabular}

\vspace{5pt}\noindent completing the result.
\end{proof}

An interesting result, which we will prove using hybrid standard form, is that when we restrict to an $\Hc$-class of $\Ren$, the partial orders $\leq^{+}$ and $\leq^{-}$ coincide. This is a decent generalization of Theorem \ref{The Partial Order Theorem} (2), as $W = H_{1}$.

\begin{prop} \label{order in a class}
For any two elements, $r, s \in \Ren$, 

(1) if $r \Lc s$ we can write their left standard forms as $r = aeb$ and $s = ceb$. Then $r \leq^{+} s$ if and only if $a \leq c$.

(2) if $r \Rc s$ we can write their right standard forms as $r = aeb$ and $s = aec$. Then $r \leq^{-} s$ if and only if $b \leq c$.

(3) if $r \Hc s$ we can write their hybrid standard forms as $r = aebec$ and $s = aedec$. Then $r \leq^{+} s$ if and only 

\indent\indent if $b \leq d$.

(4) if $r \Hc s$ we can write their hybrid standard forms as $r = aebec$ and $s = aedec$. Then $r \leq^{-} s$ if and only 

\indent\indent if $b \leq d$.
\end{prop}

\begin{proof}
(1) $r \leq^{+} s$ if and only if there exists $u \in W_{\lambda^{*}(e)}$ and $v \in W_{\lambda_{*}(e)}$ so that $a \leq cuv$ and $v^{-1}u^{-1}b \leq b$. But $v^{-1}u^{-1} \in W_{\lambda(e)}$. Since $b \in $$^{\lambda(e)}W$ it follows that $v^{-1}u^{-1}b \leq b$ if and only if $v^{-1}u^{-1} = 1$. So then $r \leq^{+} s$ if and only if $a \leq c$ as desired.

(2) is similar.

(3) Let $r = aebec$ and $s = aedec$ be the hybrid standard form decompositions and observe the following equivalences.

\vspace{5pt}\begin{tabular}{rcl}
$r \leq^{+} s$ & if and only if & $e \leq e$ and there exists $w \in W_{\lambda^{*}(e)}W_{\lambda_{*}(e)}$ so that $ab \leq adw$ and $w^{-1}c \leq c$ \\

~ & if and only if & there exists $w \in W_{\lambda(e)}$ so that $ab \leq adw$ and $w^{-1}c \leq c$\\

~ & if and only if & $ab \leq ad$ \hfill \emph{($w = 1$ since $c \leq w^{-1}c$ by Theorem \ref{Coxeter cosets} (9))} \\

~ & if and only if & $b \leq d$ \hfill \emph{(by Theorem \ref{Coxeter cosets} (3) and (6))}
\end{tabular}

\vspace{5pt}(4) is similar to (3).
\end{proof}

Interestingly, this shows that for any $r \in \Ren$ if $r \Jc e \in \Lambda$ then $(H_{r}, \leq^{+}) \cong (W_{\lambda^{*}(e)}, \leq) \cong (H_{r}, \leq^{-})$.

Later on in the paper we shall introduce the idea of minimum elements (with respect to our Adherence orders) within a given equivalence class of a given Green's relation (either $\Jc$, $\Lc$, $\Rc$, $\Hc$). However, our discussion (Section \ref{abs elts}) will take place within Renner-Coxeter monoids with a finite group of units. It turns out that some equivalences classes do have minimum elements even in the general setting. This is part of the content of our last result in this section.

\begin{thm} \label{certain abs mins exist}
Define the following three sets,

\noindent~\hfill $\GJ^{+} = \{r \in \Ren \mid r = xey\text{ in left standard form, }x = 1\}$\hfill$\JG^{-} = \{r \in \Ren \mid r = yex\text{ in right standard form, }x = 1\}$\hfill~

\noindent~\hfill$\Oh = \{r \in \Ren \mid r = aebec \text{ in hybrid standard form, } b = 1\}$\hfill~

\noindent Then,

(1) $\GJ^{+} = \{r \in \Ren \mid \forall s \in L_{r}, r \leq^{+} s \}$, $\JG^{-} = \{r \in \Ren \mid \forall s \in R_{r}, r \leq^{-} s \}$, $\Oh = \{r \in \Ren \mid \forall s \in H_{r}, r \leq^{+} s\}$, and 

\indent\indent $\Oh = \{r \in \Ren \mid \forall s \in H_{r}, r \leq^{-} s\}$.

(2) $(\GJ^{+})^{*} = \JG^{-}$ and $\Oh = \Oh^{*}$.

(3) $\GJ^{+}$, $\JG^{-}$, and $\Oh$ are monoids

(4) $\Oh$ is the smallest 

\hspace{30pt} $\bullet$ monoid containing both $\GJ^{+}$ and $\JG^{-}$

\hspace{30pt} $\bullet$ inverse monoid containing $\GJ^{+}$

\hspace{30pt} $\bullet$ inverse monoid containing $\JG^{-}$

(5) $E(\Ren) \subseteq \Oh$

(6) Suppose $r \leq^{+} s$. Then $a \leq^{+} b$ where $a \in \GJ^{+} \cap L_{r}$ and $b \in \GJ^{+} \cap L_{s}$

(7) Suppose $r \leq^{-} s$. Then $a \leq^{-} b$ where $a \in \JG^{-} \cap R_{r}$ and $b \in \JG^{-} \cap R_{s}$

(8) $\GJ^{+} \cong \Ren / \Lc$ and $\JG^{-} \cong \Ren / \Rc$

(9) $\Oh \cong \Ren / \Hc$
\end{thm}

\begin{proof}
(1) Due to the similarity between $\GJ^{+}$ and $\JG^{+}$, we shall just show $\GJ^{+} = \{r \in \Ren \mid \forall s \in L_{r}, r \leq^{+} s \}$. Suppose that $r \in \GJ^{+}$ with standard form $1ey$ by definition. Consider $s \in L_{r}$. By Proposition \ref{standard form Green's} we can see that $s = xey$ for some $x \in W^{\lambda_{*}(e)}$. Then we can see that $1 \in W_{\lambda^{*}(e)}W_{\lambda_{*}(e)}$ with $1 \leq x1$ and $1^{-1}y \leq y$, demonstrating that $r \leq^{+} s$. So $\GJ^{+} \subseteq \{r \in \Ren \mid \forall s \in L_{r}, r \leq^{+} s \}$.

Now suppose that $r \in \Ren$ and for all $s \in L_{r}$, $r \leq^{+} s$. By again invoking Proposition \ref{standard form Green's} we can write $r = xey$ and notice that $ey \Lc r$ is also in left standard form. But, by the same reasoning as the last paragraph we quickly see that $ey \leq^{+} r$. Since $r \leq^{+} ey$ by definition of $r$ we conclude that $r = ey \in \GJ^{+}$. Thus, $\GJ^{+} \supseteq \{r \in \Ren \mid \forall s \in L_{r}, r \leq^{+} s \}$ completing the result.

By the preceeding proposition, it is clear that if we show $\Oh = \{r \in \Ren \mid \forall s\in H_{r}, r \leq^{+} s\}$ then $\Oh = \{r \in \Ren \mid \forall s\in H_{r}, r \leq^{-} s\}$. For $r = ae1ec \in \Oh$ and any $s = aebec \Hc r$ we see by Proposition \ref{order in a class} $r \leq^{+} s$ if and only if $1 \leq b$ which is always true. Conversely, if $r \leq^{+} s$ for all elements in $H_{r}$ then writing them in standard form we see $r = aebec$ and $aedec$ and $b \leq d$. However, $t = ae1ec \in \Oh$ is in left standard form and clearly $t \in H_{r}$. Thus $r \leq t$, but as we have noted, $t \leq r$ and hence $r = t$. 

(2) Is clear, since $e^{*} = e$ for all idempotents, $w^{*} = w^{-1}$ for all $w \in W$ and $\Big($$^{I}W\Big)^{-1} = W^{I}$ for all $I \subseteq S$.

(3) By (2) is is clear that showing $\GJ^{+}$ is a monoid also shows that $\JG^{-}$ is a monoid. Suppose that $ex, fy \in \GJ^{+}$ are written in left standard form. Consider $exfy$ and let $x = w\circ z$ with  $w \in W^{\lambda(f)}$ (and by definition of $x$, $w \in $$^{\lambda(e)}W^{\lambda(f)}$) and $z \in W_{\lambda(f)}$. Then we find that $exfy = ewfzy = gzy$ where $g = e\wedge_{w}f$. We can see that $zy = z\circ y$ by their definitions.

Let $p = $$^{\lambda_{*}(g)}(zy)$. We claim that $p \in $$^{\lambda(g)}W$ and hence $exfy = gzy = gp \in \GJ^{+}$. Let us say $q = zyp^{-1} \in W_{\lambda_{*}(g)}$ and suppose that $sp < p$ for some $s \in \lambda^{*}(g)$. If no such $s$ exists, then $p \in $$\lambda^{*}(g)W$ and we are already done.

We can find some $t$ so that $p = st = s\circ t$. Then $zy = qp = qst = sqt$ since $q \in W_{\lambda_{*}(g)}$ and $\lambda^{*}(g)$ and $\lambda_{*}(g)$ commute. Then we see that $szy \leq zy$. Consider $sz$. There are two options, either $sz < z$, or $s > z$ (that is, $sz = s\circ z$). 

In the first case we find $z = s\tau = s\circ\tau$ for some $\tau$. It follows that $x = wz = ws\tau = sw\tau$ since $\lambda^{*}(g)$ commutes with $w \in W_{\lambda_{*)(g)}}$. This gives us a reduced word for $x$ starting with $\lambda^{*}(g) \subseteq \lambda(e)$, a contradiction. In the second case, $z < sz \in W_{\lambda^{*}(g)}W_{\lambda(f)} \subseteq W_{\lambda(f)}$ implies that $zy < (sz)y = szy$, another contradiction.

Thus, $p \in $$^{\lambda^{*}(g)}W$ and since we started with $p \in $$^{\lambda_{*}(g)}W$ we conclude $p \in $$^{\lambda^{*}(g)\sqcup\lambda_{*}(g)}W = $$^{\lambda(g)}W$ as desired. Hence $exfy = gzy = gqp = gp \in \GJ^{+}$.

To show $\Oh$ is a monoid, let $pe1eq$ be an element of $\Oh$ written in hybrid standard form. Notice that $pe1eq = (pe)(eq)$ and indeed all elements of $\Oh$ are written as $xy$ for some $x \in \GJ^{+}$ and some $y \in \JG^{-}$. To complete this result, it suffices to show that for any $x \in \GJ^{+}$ and $y \in \JG^{-}$ that $xy \in \Oh$ and $yx \in \Oh$.

Starting with $xy$, we can let $x = e\sigma$ and $y = \tau f$ for $e, f \in \Lambda$ and $\sigma \in $$^{\lambda(e)}W$, $\tau \in W^{\lambda(f)}$. Then $xy = e\sigma\tau f = \sigma(\sigma^{-1}e\sigma)(\tau f\tau^{-1})\tau = \sigma(\mu g\mu^{-1})\tau$ where $g \in \Lambda$ and $\mu \in W^{\lambda(g)}$. Since $\mu g\mu^{-1} = (\sigma^{-1} e\sigma)(\tau f\tau^{-1})$ it follows that $\mu g\mu^{-1} \leq^{+} \sigma^{-1} e\sigma, \tau f\tau^{-1}$ and $\mu g\mu^{-1} \leq^{-} \sigma^{-1} e\sigma, \tau f\tau^{-1}$ as these are all idempotents.

$\mu g \mu^{-1} \leq^{-} \sigma^{-1}e\sigma$ implies that we can find $a \in W_{\lambda^{*}(e)}$ and $\alpha \in W_{\lambda_{*}(g)}$ so that $\alpha a\sigma = \mu^{-1}$, or rather $a\sigma = \alpha^{-1}\mu^{-1}$. Similarly, we can find $b \in W_{\lambda^{*}(f)}$ and $\beta \in W_{\lambda_{*}(g)}$ so that $\tau b\beta = \mu$ or $\tau b = \mu\beta^{-1}$. Thus, $xy = \sigma\mu g\mu^{-1}\tau = a^{-1}\alpha^{-1}g\beta^{-1}b^{-1} = a^{-1}gb^{-1}$. 

Suppose that $a \not\in $$^{\lambda^{*}(g)}W$, then we can find $s \in \lambda^{*}(g)$ and a $p \in W_{\lambda^{*}(e)}$ so that $a = sp = s\circ p$ (in effect giving us a reduced word for $a$ starting with an element of $\lambda^{*}(g)$). But since $\sigma \in $$^{\lambda(e)}W$ and $\mu^{-1} \in $$^{\lambda(g)}W$ it is clear that $a\circ \sigma = a\sigma = \alpha^{-1}\mu^{-1} = \alpha^{-1}\circ \mu^{-1}$, giving us a reduced word for $\alpha^{-1}\mu^{-1}$ starting with a generator from $\lambda^{*}(g)$. However, $\alpha \in W_{\lambda_{*}(g)}$ and $\mu^{-1} \in $$^{\lambda(g)}W$ so it follows that $\alpha^{-1}\mu^{-1} \in $$^{\lambda^{*}(g)}W$ a contradiction. Thus $a^{-1} \in W^{\lambda^{*}(g)}$ and likewise $b^{-1} \in $$^{\lambda^{*}(g)}W$. It follows that $(a^{-1})^{\lambda(g)} = (a^{-1})^{\lambda_{*}(g)}$ and $^{\lambda(g)}(b^{-1}) = $$^{\lambda_{*}(g)}(b^{-1})$. So we conclude that $xy = a^{-1}gb^{-1} = (a^{-1})^{\lambda_{*}(g)}g1g$$^{\lambda_{*}(g)}()b^{-1})$ is in hybrid standard form and hence $xy \in \Oh$.

As for $yx$, we can keep our $e, f, \sigma,$ and $\tau$. Then $yx = \tau fe\sigma$. Let $g = fe = f\wedge_{1}e \in \Lambda$. Then $g \leq e, f$ and so $\lambda_{*}(e), \lambda_{*}(f) \subseteq \lambda_{*}(g)$ and $\lambda^{*}(g) \subseteq \lambda^{*}(e), \lambda^{*}(f)$. So $\lambda(g) \subseteq \lambda(e)\cup\lambda_{*}(g), \lambda(f)\cup\lambda_{*}(g)$. Thus $\tau^{\lambda_{*}(g)} \in W^{\lambda(g)}$ and $^{\lambda_{*}(g)}\sigma \in $$^{\lambda(g)}W$. So then $yx = \tau^{\lambda_{*}(g)}g$$^{\lambda_{*}(g)}\sigma = \tau^{\lambda_{*}(g)}g1g$$^{\lambda_{*}(g)}\sigma$ is in hybrid standard form, and thus $yx \in \Oh$.

So now, for any two elements $r, s \in \Oh$ we can find $x_{r}, x_{s} \in \GJ^{+}$ and $y_{r}, y_{s} \in \JG^{-}$ so that $r = x_{r}y_{r}$, $s = x_{s}y_{s}$. Then, $rs = x_{r}y_{r}x_{s}y_{s} = x_{r}x_{t}y_{t}y_{s}$ where $t = y_{r}x_{s} \in \Oh$ and $x_{t} \in \GJ^{+}$, $y_{t} \in \JG^{-}$ so that $t = x_{t}y_{t}$. Since $\GJ^{+}$ and $\JG^{-}$ are monoids we see that $rs = (x_{r}x_{t})(y_{t}y_{s})$ which is a product of something in $\GJ^{+}$ by something in $\JG^{-}$, and hence is an element of $\Oh$.

(4) Let $M$ be a monoid containing $\GJ^{+}$ and $\JG^{-}$. Let $ae1ec = aec$ be an arbitrary element of $\Oh$. Then $ae1 \in \JG^{-}$ and $1ec \in \GJ^{+}$ and so $ae1, 1ec \in M$. Since its a monoid we then see that $aec = ae1\cdot 1ec \in M$. Thus $\Oh \subseteq M$. And by (3), $\Oh$ is a monoid, so this work shows it is the smallest.

Let $M$ be an inverse monoid containing $\GJ^{+}$. Then by (2) we see that $\JG^{-} \subseteq M$ and so $\Oh \subseteq M$ by the above. Combining (2) and (3) we see that $\Oh$ is an inverse monoid, so it is indeed the smallest such monoid containing $\GJ^{+}$.

The final situation is similar to the latter.

(5) Take any idempotent, $e \in E(\Ren)$. When written in hybrid standard form we have some $f \in \Lambda$ and $x \in W^{\lambda(f)}$ so that $e = xf1fx^{-1}$. So by definition, $e \in \Oh$.

(6) Let $r = xey$ and $s = pfq$ be written in left standard form. Then $a = 1ey$ and $b = 1fq$ are written in left standard form. Since $r \leq^{+} s$ we can find $u \in W_{\lambda^{*}(f)}$ and $v \in W_{\lambda){*}(e)}$ so that $x \leq puv$ and $v^{-1}u^{-1}q \leq y$. But it is clear that for any $u, v$ we have $1 \leq 1uv$, so it follows directly that $a \leq^{+} b$.

(7) is done simiarly to result (6).

(8) As they are similar, we will only show the $\GJ^{+}$ case. Suppose that there were two elements $r, s \in \GJ^{+}$ with $r \Lc s$. Then by (1), $r \leq^{+} s \leq^{+} r$, so they are equal. Next, take any $r \in \Ren$. Write $r = xey$ in left standard form. It is clear that $r \Lc 1ey$ is also in left standard form and that $1ey \in \GJ^{+}$ by definition. Since for all $r \in \Ren$ we have $\abs{L_{r} \cap \GJ^{+}} = 1$ we conclude $\GJ^{+} \cong \Ren / \Lc$.

(9) Suppose there are two elements $r, s \in \Oh$ with $r \Hc s$. Then by (1), $r \leq^{+} s \leq^{+} r$, so then $r = s$. Next, take any $r \in \Ren$. Write $r = aebec$ in hybrid standard form. It is clear that $r \Hc ae1ec$ is also in hybrid standard form and that $ae1ec \in \Oh$. Since for all $r \in \Ren$ we have $\abs{H_{r} \cap \GJ^{+}} = 1$ we conclude $\Oh \cong \Ren / \Hc$.
\end{proof}

These monoids of minimum elements are concepts that we will return to. The existence of minimum elements only for some of Green's relations and with a particular choice of $\leq^{+}$ or $\leq^{-}$ is due to the group of units, $W$. Suppose each $\Jc$-class had a minimum element. Then for a given $e \in \Lambda$ the minimum element of $J_{e} = WeW$ must also be the minimum element of its $\Lc$-class and so has left standard form looking like $ey$. 

$ey \leq ez$ if and only if there is some $x \in W_{\lambda(e)}$ so that $xz \leq y$. But we know that $z \leq xz$, so we conclude that $z \leq y$. For an infinite $W$ this cannot be the case! With that in mind, we can explore the consequences of a finite group of units.

\section{The Finite Case and Vanilla Form} \label{vanilla form}

The more familiar object of study is Renner monoids derived from reductive algebraic monoids. Such monoids have an advantage in that they are always finite. For the remainder of this paper, any $\Ren$ considered will be finite. With a finite Renner monoid comes a finite group of units, $W$, and it is well known that a finite Coxeter group has a maximum element in its Adherence order (Theorem \ref{Coxeter length}). This element, $w_{0}$, will allow us to define a new cross-sectional lattice and with that better explore the Adherence order and its relationship with Green's relations.

\begin{prop}
Define the set of idempotents, $\Lambda^{-} := w_{0}\Lambda w_{0}$, called the {\bf opposite cross-sectional lattice}. Then $(\Ren, \Lambda^{-}, S)$ is a Renner-Coxeter system.
\end{prop}

\begin{proof}
Properties, (i) and (iii) are not altered by $\Lambda^{-}$.

(ii) Let $e \in E(\Ren)$. Since $\Lambda$ is a transversal we can find $u \in W$ and $f \in \Lambda$ so that $e = ufu^{-1}$. But then, $w_{0}fw_{0} \in \Lambda^{-}$ by definition and it follows that $e = (uw_{0})(w_{0}fw_{0})(w_{0}u^{-1}) = (uw_{0})(w_{0}fw_{0})(uw_{0})^{-1}$. Furthermore, if $e \Jc f$ for two $e, f \in \Lambda^{-}$ then we can find $w \in W$ so that $wew^{-1} = f$ and $g, h \in \Lambda$ so that $w_{0}ew_{0} = g$ and $w_{0}fw_{0} = h$. Combining these we see that $w_{0}ww_{0}gw_{0}w^{-1}w_{0} = h$. But $\Lambda$ is a transversal, so $g = h$ and hence $e = f$.

(iv) Suppose that $e \leq e^{\prime} \in E(\Ren)$. Then $ee^{\prime} = e^{\prime} = e^{\prime}e$. Since $(\Ren, \Lambda, S)$ is a Renner-Coxeter system we can find $f \leq f^{\prime} \in \Lambda$ and $w \in W$ so that $wew^{-1} = f$ and $we^{\prime}w^{-1} = f^{\prime}$. 

But then $(w_{0}w)e(w_{0}w)^{-1} = w_{0}wew^{-1}w_{0} = w_{0}fw_{0} = g \in \Lambda^{-}$ and $(w_{0}w)e^{\prime}(w_{0}w)^{-1} = w_{0}we^{\prime}w^{-1}w_{0} = w_{0}f^{\prime}w_{0} = g^{\prime} \in \Lambda^{-}$. We can also notice that $f \leq f^{\prime}$ implies $ff^{\prime} = f = f^{\prime}f$ and so $gg^{\prime} = w_{0}ff^{\prime}w_{0} = w_{0}fw_{0}w_{0}f^{\prime}w_{0} = w_{0}fw_{0} = g = w_{0}fw_{0} = w_{0}f^{\prime}fw_{0} = w_{0}f^{\prime}w_{0}w_{0}fw_{0} = g^{\prime}g$. So $g \leq g^{\prime} \in \Lambda^{-}$ and $(w_{0}w) \in W$ satisfy this requirement.

(v) Let $e \in \Lambda^{-}$ and $f = w_{0}ew_{0} \in \Lambda$. Notice that because $\ell(w_{0}ww_{0}) = \ell(w)$ and hence, $w_{0}Sw_{0} = S$ we have the series of equalities, $\lambda^{*}(e) = \{s \in S \mid se = es \neq e\} = \{s \in S \mid w_{0}sw_{0}w_{0}ew_{0} = w_{0}ew_{0}w_{0}sw_{0} \neq w_{0}ew_{0}\} = w_{0}\{s \in S \mid s(w_{0}ew_{0}) = (w_{0}ew_{0})s \neq w_{0}ew_{0}\}w_{0} = w_{0}\lambda^{*}(w_{0}ew_{0})w_{0} \subseteq S$. So the map $\lambda^{*} : \Lambda^{-} \rightarrow P(S)$ is certainly defined.

Pick any $e \leq e^{\prime} \in \Lambda$ and let $f = w_{0}ew_{0} \in \Lambda$, $f^{\prime} = w_{0}e^{\prime}w_{0} \in \Lambda$. As we noted in (iv) $e \leq e^{\prime}$ implies that $f \leq f^{\prime}$. So then $\lambda^{*}(e) = w_{0}\lambda^{*}(f)w_{0} \subseteq w_{0}\lambda^{*}(f^{\prime})w_{0} = \lambda^{*}(e^{\prime})$ as desired.

(vi) Notice, as we similarly saw in (v), that, using the fact of $w_{0}Sw_{0} = S$, for any $e \in \Lambda^{-}$, we have the chain $\{w \in W \mid we = ew\} = w_{0}\{w \in W \mid w_{0}ww_{0}w_{0}ew_{0} = w_{0}ew_{0}w_{0}ww_{0}\}w_{0} = w_{0}W_{\lambda(w_{0}ew_{0})}w_{0} = W_{w_{0}\lambda(w_{0}ew_{0})w_{0}}$ and  \scalebox{0.95}{$\{w \in W \mid we = ew = e\}$ $= w_{0}\{w \in W \mid w_{0}ww_{0}w_{0}ew_{0} = w_{0}ew_{0}w_{0}ww_{0} = w_{0}ew_{0}\}w_{0} = w_{0}W_{\lambda(w_{0}ew_{0})}w_{0} = W_{w_{0}\lambda_{*}(w_{0}ew_{0})w_{0}}$}. So both are clearly standard parabolic subgroups of $W$.
\end{proof}

This means that we know have four different standard forms for a given element of $\Ren$ and four different partial orders. Our next result relates partial orders in $(\Ren, \Lambda, S)$ and $(\Ren, \Lambda^{-}, S)$.

\begin{prop} \label{order and w_0}
For any two elements, $r, s\in \Ren$, the following are equivalent.

(1) $r \leq^{+} s$ in $(\Ren, \Lambda, S)$

(2) $w_{0}rw_{0} \leq^{+} w_{0}sw_{0}$ in $(\Ren, \Lambda^{-}, S)$
\end{prop}

\begin{proof}
Let $r = xey$ and $s = afb$ be written in left standard form. We claim $w_{0}rw_{0} = (w_{0}xw_{0})(w_{0}ew_{0})(w_{0}yw_{0})$ and $w_{0}rw_{0} = (w_{0}aw_{0})(w_{0}fw_{0})(w_{0}bw_{0})$ are in left standard form in $(\Ren, \Lambda^{-}, S)$. Suppose $(w_{0}xw_{0}) \not\in W^{\lambda_{*}(w_{0}ew_{0})}$, then we can find a reduced word expression $(w_{0}xw_{0}) = s_{1}\cdots s_{k}$ with $s_{k} \in \lambda_{*}(w_{0}ew_{0})$. However, for any $s \in S$, $\ell(w_{0}sw_{0}) = \ell(s) = 1$, so $w_{0}sw_{0}$ is a generator and it follows that $(w_{0}s_{1}w_{0})\cdots (w_{0}s_{k}w_{0})$ is a reduced word for $x$ ending in $(w_{0}s_{k}w_{0}) \in w_{0}\lambda_{*}(w_{0}ew_{0})w_{0} = \lambda_{*}(e)$, a contradiction. Similar arguments for $(w_{0}yw_{0})$, $(w_{0}aw_{0})$, and $(w_{0}bw_{0})$ establish our claim.

Then we note the following equivalences, 

\begin{tabular}{rcl}
$r \leq^{+} s$ & if and only if & $e \leq f$ and there exists $w \in W_{\lambda^{*}(f)}W_{\lambda_{*}(e)}$ so that $x \leq aw$ and $w^{-1}b \leq y$ \\

~ & if and only if & $w_{0}ew_{0} \leq w_{0}fw_{0}$ and there exists $w \in w_{0}W_{\lambda^{*}(w_{0}fw_{0})}w_{0}w_{0}W_{\lambda_{*}(w_{0}ew_{0})}w_{0}$ so \\

~ & ~ & \indent\indent\indent that $w_{0}xw_{0} \leq w_{0}aw_{0}w_{0}ww_{0}$ and $w_{0}w^{-1}w_{0}w_{0}bw_{0} \leq w_{0}yw_{0}$ \\

~ & if and only if & $w_{0}ew_{0} \leq w_{0}fw_{0}$ and there exists $w \in W_{\lambda^{*}(w_{0}fw_{0})}W_{\lambda_{*}(w_{0}ew_{0})}$ so \\

~ & ~ & \indent\indent\indent that $w_{0}xw_{0} \leq w_{0}aw_{0}w$ and $w^{-1}w_{0}bw_{0} \leq w_{0}yw_{0}$ \\

~ & if and only if & $w_{0}rw_{0} \leq^{+} w_{0}sw_{0}$
\end{tabular}

\noindent completing the result.
\end{proof}

\begin{thm}
Let $r \in \Ren$. There exists a unique set of elements, $\sigma_{-}, e_{-}, \sigma_{0}, e_{+}, \sigma_{+} \in \Ren$ satisfying the following four properties:

(i) $r = \sigma_{-}e_{-}\sigma_{0}e_{+}\sigma_{+}$

(ii) $e_{+} \in \Lambda$, $e_{-} \in \Lambda^{-}$, and $e_{-} \Jc r \Jc e_{+}$ 

(iii) $\sigma_{+} \in~$$^{\lambda(e_{+})}W$ and $\sigma_{-} \in W^{\lambda(e_{-})}$

(iv) $\sigma_{0} \in W_{\lambda^{*}(e_{-})}\Big($$^{\lambda_{*}(e_{-})}w_{0}^{\lambda_{*}(e_{+})}\Big)W_{\lambda^{*}(e_{+})} = $$^{\lambda_{*}(e_{-})}\Big(W_{\lambda(e_{-})}w_{0}W_{\lambda(e_{+})}\Big)^{\lambda_{*}(e_{+})}$

\noindent When we write the element $r = \sigma_{-}e_{-}\sigma_{0}e_{+}\sigma_{+}$ in this way we say that $r$ is written in {\bf vanilla form}.
\end{thm}

\begin{proof}
Let $r = xey$ be in left standard form with respect to $(\Ren, \Lambda, S)$. Notice that $f = w_{0}ew_{0} \in \Lambda^{-}$ and so $r = xey = xeey = xw_{0}fw_{0}ey$. Let $z = (xw_{0})^{\lambda(f)}$ then we can find $u \in W_{\lambda(f)}$ so that $zu = xw_{0}$ and by substituting in we get, $r = zufw_{0}ey = zfuw_{0}ey$. Now, $uw_{0} \in W_{\lambda(f)}w_{0} = W_{\lambda(f)}W_{\lambda(f)}w_{0} = W_{\lambda(f)}w_{0}W_{\lambda(e)}w_{0}w_{0} = W_{\lambda(f)}w_{0}W_{\lambda(e)}$. By taking $v = $$^{\lambda_{*}(f)}(uw_{0})^{\lambda_{*}(e)}$ we see that $r = zfvey$ satisfies the conditions of our vanilla form, establishing existence.

Toward uniqueness, suppose that $r = \sigma_{-}e_{-}\sigma_{0}e_{+}\sigma_{+} = \tau_{-}f_{-}\tau_{0}f_{+}\tau_{+}$ and both decompositions satisfy the vanilla form conditions. Then $e_{-} \Jc f_{-}$ and $e_{+} \Jc f_{+}$, and since $\Lambda^{-}$ and $\Lambda$ are transversals it is immediate that $e_{-} = f_{-}$ and $e_{+} = f_{+}$. 

Since $\sigma_{0}, \tau_{0} \in $$^{\lambda_{*}(e_{-})}\Big(W_{\lambda(e_{-})}w_{0}W_{\lambda(e_{+})}\Big)^{\lambda_{*}(e_{+})} \subseteq W_{\lambda(e_{-})}w_{0}W_{\lambda(e_{+})}$ we can quickly see that $e_{-}\sigma_{0} = e_{-}\sigma_{0}e_{+} = \sigma_{0}e_{+}$ and $e_{-}\tau_{0} = e_{-}\tau_{0}e_{+} = \tau_{0}e_{+}$. Then $(\sigma_{-}\sigma_{0})^{\lambda_{*}(e_{+})}e_{+}\sigma_{+} = (\sigma_{-}\sigma_{0})e_{+}\sigma_{+} = r = (\tau_{-}\tau_{0})e_{+}\tau_{+} = (\tau_{-}\tau_{0})^{\lambda_{*}(e_{+})}e_{+}\tau_{+}$ and the ends are in left standard form. By uniqueness of left standard form then $\sigma_{+} = \tau_{+}$. Similarly, $\sigma_{-} = \tau_{-}$.

So then $\sigma_{-}^{-1}r\sigma_{+}^{-1} = e_{-}\sigma_{0}e_{+} = e_{-}\tau_{0}e_{+}$, or rather $e_{-}\sigma_{0} = e_{-}\tau_{0}$. Rearranging we find that $e_{-}\sigma_{0}\tau_{0}^{-1} = e_{-}$ and hence $\sigma_{0}\tau_{0}^{-1} \in W_{\lambda_{*}(e_{-})}$. So then $\sigma_{0} \in W_{\lambda_{*}(e_{-})}\tau_{0} \subseteq W_{\lambda_{*}(e_{-})}\tau_{0}W_{\lambda_{*}(e_{+})}$. But since $\sigma_{0}, \tau_{0} \in $$^{\lambda_{*}(e_{-})}W^{\lambda_{*}(e_{+})}$ we conclude that $\sigma_{0} = \tau_{0}$.
\end{proof}

Vanilla form (named to avoid overused terms like normal and canonical) can be interestingly contrasted with hybrid standard form. Both bridge the gap between left and right standard forms. Both (as we will see in a bit) allow us to identify each of Green's relations by observing the decomposition. However, the hybrid standard form is contained solely in $(\Ren, \Lambda, S)$, whereas our new vanilla form allows us to study both $(\Ren, \Lambda, S)$ and $(\Ren, \Lambda^{-}, S)$.

\begin{cor} \label{standard from vanilla}
Suppose $r = \sigma_{-}e_{-}\sigma_{0}e_{+}\sigma_{+}$ is in vanilla form. Then,

(1) $\sigma_{+}^{-1} \leq_{L} \sigma_{0}$ and $\sigma_{-}^{-1} \leq_{R} \sigma_{0}$

(2) $r = (\sigma_{-}\sigma_{0})e_{+}\sigma_{+}$ is in left standard form (in $(\Ren, \Lambda, S)$) 

(3) $r = \sigma_{-}e_{-}(\sigma_{0}\sigma_{+})$ is in right standard form (in $(\Ren, \Lambda^{-}, S)$)
\end{cor}

\begin{proof}
(1) We shall just show $\sigma_{+}^{-1} \leq_{L} \sigma_{0}$ as $\sigma_{-}^{-1} \leq_{R} \sigma_{0}$ is done symmetrically. $\sigma_{+} \in $$^{\lambda(e_{+})}W$, so $\sigma_{+}^{-1} \in W^{\lambda(e_{+})}$. $\sigma_{+}^{-1}w_{0}(\lambda(e_{+})) = \sigma_{+}^{-1}\circ w_{0}(\lambda(e_{+})) \leq_{L} w_{0}$ (as $w \leq_{L} w_{0}$ for all $w \in W$). This means we can find $u \in W$ so that $u\sigma_{+}^{-1}w_{0}(\lambda(e_{+})) = w_{0}$ with $\ell(u) + \ell(\sigma_{+}^{-1}w_{0}(\lambda(e_{+}))) = \ell(w_{0})$. But then, since $\ell(\sigma_{+}^{-1}w_{0}(\lambda(e_{+}))) = \ell(\sigma_{+}^{-1}) + \ell(w_{0}(\lambda(e_{+})))$ and $w_{0} = w_{0}^{\lambda(e_{+})}\circ w_{0}(\lambda(e_{+})) = w_{0}^{\lambda(e_{+})}w_{0}(\lambda(e_{+}))$ it follows that $w_{0}^{\lambda(e_{+})} = u\sigma_{+}^{-1}$ and $\ell(w_{0}^{\lambda(e_{+})}) = \ell(u) + \ell(\sigma_{+}^{-1})$. Thus $\sigma_{+}^{-1} \leq_{L} w_{0}^{\lambda(e_{+})}$.

We know that $\sigma_{0} \in W_{\lambda(e_{-})}w_{0}W_{\lambda(e_{+})} = W_{\lambda(e_{-})}$$^{\lambda(e_{-})}w_{0}^{\lambda(e_{+})}W_{\lambda(e_{+})} = W_{\lambda(e_{-})}$$^{\lambda(e_{-})}w_{0}^{\lambda(e_{+})}$. Now, since $^{\lambda(e_{-})}w_{0}^{\lambda(e_{+})} = w_{0}^{\lambda(e_{+})}$ it follows that we can find $v \in W_{\lambda(e_{-})}$ so that $\sigma_{0} = vw_{0}^{\lambda(e_{+})}$ and $\ell(\sigma_{0}) = \ell(v) + \ell(w_{0}^{\lambda(e_{+})})$. So $w_{0}^{\lambda(e_{+})} \leq_{L} \sigma_{0}$ and since $\leq_{L}$ is transitive, $\sigma_{+}^{-1} \leq_{L} w_{0}^{\lambda(e_{+})} \leq_{L} \sigma_{0}$ implies $\sigma_{+}^{-1} \leq_{L} \sigma_{0}$.

(2) and (3) both follow from (1). Since $\sigma_{-}^{-1} \leq_{R} \sigma_{0}$ we see that $\sigma_{0} = \sigma_{-}^{-1}w = \sigma_{-}^{-1}\circ w$ for some $w \in W$. Thus, $r = \sigma_{-}\sigma_{0}e_{+}\sigma_{+} = we\sigma_{+}$, so we only need show that $w \in W^{\lambda_{*}(e)}$. Suppose not, then there exists a reduced word for $w$ which ends in an element of $\lambda_{*}(e)$. However, then $\sigma_{-}^{-1}w$ gives a reduced word expression for $\sigma_{0}$ ending in the same element of $\lambda_{*}(e)$, a contradiction. This demonstrates (2) and (3) is similar.
\end{proof}

\begin{cor} \label{Green's relations and vanilla form}
Suppose $r = \sigma_{-}e_{-}\sigma_{0}e_{+}\sigma_{+}$ and $r = \tau_{-}f_{-}\tau_{0}f_{+}\tau_{+}$ are in vanilla form. Then,

(1) $r\Jc s$ if and only if $e_{-} = f_{-}$ if and only if $e_{+} = f_{+}$

(2) $r\Lc s$ if and only if $e_{-} = f_{-}$, $e_{+} = f_{+}$, and $\sigma_{+} = \tau_{+}$

(3) $r\Rc s$ if and only if $\sigma_{-} = \tau_{-}$, $e_{-} = f_{-}$, and $e_{+} = f_{+}$

(4) $r\Hc s$ if and only if $\sigma_{-} = \tau_{-}$, $e_{-} = f_{-}$, $e_{+} = f_{+}$, and $\sigma_{+} = \tau_{+}$
\end{cor}

\begin{proof}
(1) We know we can write $r = (\sigma_{-}\sigma_{0})e_{+}\sigma_{+}$ and $s = (\tau_{-}\tau_{0})f_{+}\tau_{+}$ in left standard form. By Proposition \ref{standard form Green's} we know $r \Jc s$ if and only if $e_{+} = f_{+}$. Similarly, we can show $r \Jc s$ if and only if $e_{-} = f_{-}$ using right standard form.

(2) We can write $r = (\sigma_{-}\sigma_{0})e_{+}\sigma_{+}$ and $s = (\tau_{-}\tau_{0})f_{+}\tau_{+}$ in left standard form. By Proposition \ref{standard form Green's} we know $r \Lc s$ if and only if $e_{+} = f_{+}$ and $\sigma_{+} = \tau_{+}$.

(3) is done similar to (2).

(4) is a combination of (2) and (3).
\end{proof}

The most important result of vanilla form is the following, allowing us to show an equivalence between Adherence orders on $(\Ren, \Lambda, S)$ and $(\Ren, \Lambda^{-}, S)$.

\begin{thm} \label{vanilla order}
Suppose $r = \sigma_{-}e_{-}\sigma_{0}e_{+}\sigma_{+}$ and $r = \tau_{-}f_{-}\tau_{0}f_{+}\tau_{+}$ are in vanilla form. Then the following are equivalent,

(1) $e_{-} \leq f_{-}$, $e_{+} \leq f_{+}$, and there exist $w_{-} \in W_{\lambda_{*}(e_{-})}W_{\lambda^{*}(f_{-})}$ and $w_{+} \in W_{\lambda^{*}(f_{+})}W_{\lambda_{*}(e_{+})}$ so that $\tau_{-}w_{-}^{-1} \leq \sigma_{-}$, $\sigma_{0} \leq w_{-}\tau_{0}w_{+}$, and $w_{+}^{-1}\tau_{+} \leq \sigma_{+}$

(2) $r \leq^{+} s$ as elements of $(\Ren, \Lambda, S)$

(3) $r \leq^{-} s$ as elements of $(\Ren, \Lambda^{-}, S)$
\end{thm}

\begin{proof}
(2) $\Rightarrow$ (1) Suppose we have $w_{+} \in W_{\lambda^{*}(f_{+})}W_{\lambda_{*}(e_{+})}$ such that $\sigma_{-}\sigma_{0} \leq \tau_{-}\tau_{0}w_{+}$, and $w_{+}^{-1}\tau_{+} \leq \sigma_{+}$. We can see that $\tau_{0}w_{+} \in W_{\lambda(f_{-})}w_{0}W_{\lambda(f_{+})} = W_{\lambda(f_{-})}w_{0}$, so there exists $y \in W_{\lambda(f_{-})}$ such that $\tau_{0}w_{+} = yw_{0}$. Likewise there exists $x \in W_{\lambda(e_{-})}$ so that $\sigma_{0} = xw_{0}$. Thus, $\sigma_{-}xw_{0} = \sigma_{-}\sigma_{0} \leq \tau_{-}\tau_{0}w_{+} = \tau_{-}yw_{0}$ which is equivalent to $\tau_{-}y \leq \sigma_{-}x$.

Since $\sigma_{-} \in W^{\lambda(e_{-})}$ and $x \in W_{\lambda(e_{-})}$, by Theorem \ref{Coxeter cosets} we can find $a, b, c \in W$ so that $a = (\tau_{-}y)^{\lambda(e_{-})}$, $b, c \in W_{\lambda(e_{-})}$ with $abc = a\circ b\circ c$, $ab \leq \sigma_{-}$ and $c \leq x$. Let $w_{-} = cy^{-1} \in W_{\lambda(e_{-})}W_{\lambda(f_{-})} = W_{\lambda_{*}(e_{-})}W_{\lambda^{*}(f_{-})}$. Then we see that $\tau_{-}w_{-}^{-1} = \tau_{-}yc^{-1} = abcc^{-1} = ab \leq \sigma_{-}$ and $w_{-}y = c \leq x$ $\Leftrightarrow$ $\sigma_{0} = xw_{0} \leq w_{-}yw_{0} = w_{-}\tau_{0}w_{+}$.

(1) $\Rightarrow$ (2) Suppose we have $w_{-} \in W_{\lambda_{*}(e_{-})}W_{\lambda^{*}(f_{-})}$ and $w_{+} \in W_{\lambda^{*}(f_{+})}W_{\lambda_{*}(e_{+})}$ such that $\tau_{-}w_{-}^{-1} \leq \sigma_{-}$, $\sigma_{0} \leq w_{-}\tau_{0}w_{+}$, and $w_{+}^{-1}\tau_{+} \leq \sigma_{+}$. By Theorem \ref{Coxeter -fixes}, since $\sigma_{-}^{-1} \leq_{R} \sigma_{0}$, we can conclude that $\sigma_{-}\sigma_{0} \leq \tau_{-}w_{-}^{-1}w_{-}\tau_{0}w_{+} = \tau_{-}\tau_{0}w_{+}$ which satisfies the conditions of $r \leq^{+} s$.

(3) $\Leftrightarrow$ (1) follows by nearly identical reasoning to our previous two arguments, finishing off this result.
\end{proof}

The following corollary gives us a set of broader results than those in Proposition \ref{order in a class}.

\begin{cor} \label{Green's relations and the Adherence order}
Take two elements $r, s \in \Ren$ and let their vanilla forms be $r = \sigma_{-}e_{-}\sigma_{0}e_{+}\sigma_{+}$ and $s = \tau_{-}f_{-}\tau_{0}f_{+}\tau_{+}$.

(1) If $r \Jc s$ then $r\leq^{+} s$ if and only if there exist elements $w_{-} \in W_{\lambda^{*}(e_{-})}$ and $w_{+} \in W_{\lambda^{*}(e_{+})}$ so that $\tau_{-}w_{-}^{-1} \leq \sigma_{-}$, $\sigma_{0} \leq w_{-}\tau_{0}w_{+}$, and $w_{+}^{-1}\tau_{+} \leq \sigma_{+}$

(2) If $r \Lc s$ then $r\leq^{+} s$ if and only if $\sigma_{-}\sigma_{0} \leq \tau_{-}\tau_{0}$

(3) If $r \Rc s$ then $r\leq^{+} s$ if and only if $\sigma_{0}\sigma_{+} \leq \tau_{0}\tau_{+}$

(4) If $r \Hc s$ then $r\leq^{+} s$ if and only if $\sigma_{0} \leq \tau_{0}$
\end{cor}

\begin{proof}
(1) Since $r \Jc s$ Corollary \ref{Green's relations and vanilla form} tells us that $e_{-} = f_{-}$ and $e_{+} = f_{+}$, and so our condition for $r \leq^{+} s$ reduces to $w_{-} \in W_{\lambda(e_{-})}$, $w_{+} \in W_{\lambda(e_{+})}$ so that $\tau_{-}w_{-}^{-1} \leq \sigma_{-}$, $\sigma_{0} \leq w_{-}\tau_{0}w_{+}$, $w_{+}^{-1}\tau_{+} \leq \sigma_{+}$. Let $u_{+} \in W_{\lambda^{*}(e_{+})}$, $v_{+} \in W_{\lambda_{*}(e_{+})}$, $u_{-} \in W_{\lambda^{*}(e_{-})}$, $v_{-} \in W_{\lambda_{*}(e_{-})}$ so that $w_{+} = u_{+}v_{+} = u_{+}\circ v_{+}$ and $w_{-} = v_{-}u_{-} = v_{-}\circ u_{-}$. We can quickly see that $\tau_{-}u_{-}^{-1} \leq \tau_{-}w_{-}^{-1} \leq \sigma_{-}$, $\sigma_{0} \leq u_{-}\tau_{0}u_{+}$, and $u_{+}^{-1}\tau_{+} \leq w_{+}^{-1}\tau_{+} \leq \sigma_{+}$.

The reverse implication is obvious.

(2) Since $r\Lc s$ then $e_{+} = f_{+}$ and $\sigma_{+} = \tau_{+}$. By Corollary \ref{standard from vanilla} we know that $r = (\sigma_{-}\sigma_{0})e_{+}\sigma_{+}$ and $s = (\tau_{-}\tau_{0})e_{+}\sigma_{+}$ are in left standard form. Then $r \leq^{+} s$ if and only if we can find $w_{+} \leq W_{\lambda^{*}(e_{+})}W_{\lambda_{*}(e_{+})}$ so that $(\sigma_{-}\sigma_{0}) \leq (\tau_{-}\tau_{0})w_{+}$ and $w_{+}^{-1}\sigma_{+} \leq \sigma_{+}$. Now, since $w_{+} \in W_{\lambda(e_{+})}$ we see that $\sigma_{+} \leq w_{+}^{-1}\sigma_{+} = \sigma_{+}$ forcing $w_{+} = 1$. Thus, $r \leq^{+} s$ if and only if $\sigma_{-}\sigma_{0} \leq \tau_{-}\tau_{0}$.

(3) is similar to (2).

(4) Since $r \Hc s$ we know by Corollary \ref{Green's relations and vanilla form} that $\sigma_{-} = \tau_{-}$, $e_{-} = f_{-}$, $e_{+} = f_{+}$, and $\sigma_{+} = \tau_{+}$. $r \leq^{+} s$ if and only if there exists $w_{-} \in W_{\lambda(e_{-})}$ and $w_{+} \in W_{\lambda(e_{+})}$ so that $\sigma_{-}w_{-}^{-1} \leq \sigma_{-}$, $\sigma_{0} \leq w_{-}\tau_{0}w_{+}$, and $w_{+}^{-1}\sigma_{+} \leq \sigma_{+}$. Like (2) it is quick to show $w_{-} = w_{+} = 1$, reducing our condition to $r \leq^{+} s$ if and only if $\sigma_{0} \leq \tau_{0}$.
\end{proof}

Of course, one could construct a nearly identical decomposition to the vanilla form which allows for the quick comparisons of Green's relations but instead allows the computation of $\leq^{-}$ and extends right standard form. This, ``{\bf chocolate}'' form while an interesting concept perhaps, would be functionally no different from vanilla form (indeed one can derive this new form by taking the inverse of the vanilla form for $r^{*}$, or simply studying vanilla form in $(\Ren, \Lambda^{-}, S)$), and so we will confine our results to focussing on the vanilla form. This will mean placing emphasis on $\leq^{+}$ over $\leq^{-}$, but interested readers can easily mimic the proofs and derive the corresponding results.

\begin{cor} \label{Adherence equivalences}
The following are equivalent for any $r, s \in \Ren$

(1) $r \leq^{+} s$ in $(\Ren, \Lambda, S)$

(2) $r^{*} \leq^{-} s^{*}$ in $(\Ren, \Lambda, S)$

(3) $w_{0}rw_{0} \leq^{+} w_{0}sw_{0}$ in $(\Ren, \Lambda^{-}, S)$

(4) $r \leq^{-} s$ in $(\Ren, \Lambda^{-}, S)$

(5) $r^{*} \leq^{+} s^{*}$ in $(\Ren, \Lambda^{-}, S)$

(6) $w_{0}rw_{0} \leq^{-} w_{0}sw_{0}$ in $(\Ren, \Lambda, S)$
\end{cor}

This corollary merely sums up the various Adherence order equivalences we have come across.

\begin{proof}
(1) $\Leftrightarrow$ (2) and (4) $\Leftrightarrow$ (5) follow from Proposition \ref{order and *}. (1) $\Leftrightarrow$ (3) and (4) $\Leftrightarrow$ (6) follow from Proposition \ref{order and w_0}. And lastly, (1) $\Leftrightarrow$ (4) by Theorem \ref{vanilla order}.
\end{proof}

We will close out this section with a very interesting result about the $\Jc$ relation and the Adherence order. Recalling that $r \Jc s$ if and only if there exists $t$ so that $r \Rc t\Lc s$ if and only if there exists $u$ so that $r \Lc u \Rc s$, the following proposition tells us that these $t$ and $u$ may be chosen with the Adherence order in mind.

\begin{prop}
For any two elements, $r, s \in \Ren$ with $r\Jc s$ and $r \leq^{+} s$, there exist elements, $t, u \in \Ren$ so that $r \leq^{+} t \leq^{+} s$, $r\leq^{+} u \leq^{+} s$ and $r \Rc t \Lc s$ and $r \Lc u \Rc s$.
\end{prop}

\begin{proof}
Let $r = \sigma_{-}e_{-}\sigma_{0}e_{+}\sigma_{+}$ and $s = \tau_{-}e_{-}\tau_{0}e_{+}\tau_{+}$ be in vanilla form. By Corollary \ref{Green's relations and the Adherence order} since $r \leq^{+} s$ there exist $w_{-} \in W_{\lambda^{*}(e_{-})}$ and $w_{+} \in W_{\lambda^{*}(e_{+})}$ such that $\tau_{-}w_{-}^{-1} \leq \sigma_{-}$, $\sigma_{0} \leq w_{-}\tau_{0}w_{+}$, and $w_{+}^{-1}\tau_{+} \leq \sigma_{+}$. Let $t = \sigma_{-}e_{-}(w_{-}\tau_{0})e_{+}\tau_{+}$ and $u = \tau_{-}e_{-}(\tau_{0}w_{+})e_{+}\sigma_{+}$. It is clear from the definition of vanilla form that $\tau_{0} \in W_{\lambda^{*}(e_{-})}\Big($$^{\lambda_{*}(e_{-})}w_{0}^{\lambda_{*}(e_{+})}\Big)W_{\lambda^{*}(e_{+})}$ and so $(w_{-}\tau_{0}), (\tau_{0}w_{+}) \in W_{\lambda^{*}(e_{-})}\Big($$^{\lambda_{*}(e_{-})}w_{0}^{\lambda_{*}(e_{+})}\Big)W_{\lambda^{*}(e_{+})}$, hence $t$ and $u$ are in vanilla form.

Since they are all in vanilla form, Corollary \ref{Green's relations and vanilla form} allows us to easily conclude that $r \Rc t \Lc s$ and $r \Lc u \Rc s$. Observe that $1 \in W_{\lambda^{*}(e_{-})}$ and $w_{+} \in W_{\lambda^{*}(e_{+})}$ satisfy $\sigma_{-}1^{-1} \leq \sigma_{-}$, $\sigma_{0} \leq 1(w_{-}\tau_{0})w_{+}$, and $w_{+}^{-1}\sigma_{+} \leq \tau_{+}$. So $r \leq^{+} t$. $w_{-}^{-1} \in W_{\lambda^{*}(e_{-})}$ and $1 \in W_{\lambda^{*}(e_{+})}$ satisfy $\tau_{-}w_{-}^{-1} \leq \sigma_{-}$, $(w_{-}\tau_{0}) \leq w_{-}\tau_{0}1$, and $1^{-1}\tau_{+} \leq \tau_{+}$. So $t \leq^{+} s$. Likewise we can show $r \leq^{+} u \leq^{+} s$, completing the proof.
\end{proof}

\section{Maximum and Minimum Elements} \label{abs elts}

In Section 5 of \cite{PPR} an interesting element of a given $\Hc$-class in a Renner monoid is distinguished. For  This particular example motivates the following definition and this section's discussion.

\begin{defn}
For an equivalence relation on $\Ren$, $\Tc$, $\varepsilon = +$ or $-$, and any element $r \in \Ren$ we define,

\noindent~\hfill $max^{\varepsilon}T_{r} = \left\{\begin{array}{cl}s \in T_{r}&\text{if }\forall t\in T_{r}, t\leq^{\varepsilon} s\\\text{undefined} &\text{otherwise}\end{array}\right.$ \hfill $min^{\varepsilon}T_{r} = \left\{\begin{array}{cl}s \in T_{r}&\text{if }\forall t\in T_{r}, s\leq^{\varepsilon} t\\\text{undefined} &\text{otherwise}\end{array}\right.$ \hfill~

The left element, if it exists, is called the {\bf absolute maximum element} of the $\Tc$-class of $r$. The right element, if it exists, is called the {\bf absolute minimum element} of the $\Tc$-class.
\end{defn}

\begin{prop}
Suppose there are two equivalence relations, $\Tc$ and $\Uc$ so that $r \Tc s$ implies $r \Uc s$ for all $r, s \in \Ren$. Then for any $r \in \Ren$, assuming they exist, $max^{\varepsilon}T_{r} \leq^{\varepsilon} max^{\varepsilon}U_{r}$ and $min^{\varepsilon}U_{r} \leq^{\varepsilon} min^{\varepsilon}T_{r}$.
\end{prop}

\begin{proof}
Suppose that $t = max^{\varepsilon}T_{r}$ and $u = max^{\varepsilon}U_{r}$. Then it is clear that $t \Tc r \Uc u$ implies that $t \Uc u$. But by definition $s \leq^{\varepsilon} u$ for all $s \Uc u$. Then $t \leq^{\varepsilon} u$ as desired. The minimum case is done similarly.
\end{proof}

As it turns out, for the Green's relations, $\Hc$, $\Lc$, $\Rc$, and $\Jc$ (in a finite $\Ren$) these maximum and minimum elements always exist and can be characterized using the vanilla form decomposition.

\begin{thm} \label{abs elts vanilla conditions}
Let $r = \sigma_{-}e_{-}\sigma_{0}e_{+}\sigma_{+}$ be written in vanilla form in $(\Ren, \Lambda^{\varepsilon}, S)$.

(1) $r = min^{\epsilon}H_{r}$ if and only if $\sigma_{0} = $$^{\lambda(e_{-})}w_{0}\phantom{~\hspace{-1pt}}^{\lambda(e_{+})}$

(2) $r = min^{\epsilon}L_{r}$ if and only if $\sigma_{-} = $$^{\lambda(e_{+})}w_{0}\phantom{~\hspace{-1pt}}^{\lambda(e_{-})}$, $\sigma_{0} = $$^{\lambda(e_{-})}w_{0}\phantom{~\hspace{-1pt}}^{\lambda(e_{+})}$

(3) $r = min^{\epsilon}R_{r}$ if and only if $\sigma_{0} = $$^{\lambda(e_{-})}w_{0}\phantom{~\hspace{-1pt}}^{\lambda(e_{+})}$, $\sigma_{+} = $$^{\lambda(e_{+})}w_{0}\phantom{~\hspace{-1pt}}^{\lambda(e_{-})}$

(4) $r = min^{\epsilon}J_{r}$ if and only if $\sigma_{-} = $$^{\lambda(e_{+})}w_{0}\phantom{~\hspace{-1pt}}^{\lambda(e_{-})}$, $\sigma_{0} = $$^{\lambda(e_{-})}w_{0}\phantom{~\hspace{-1pt}}^{\lambda(e_{+})}$, $\sigma_{+} = $$^{\lambda(e_{+})}w_{0}\phantom{~\hspace{-1pt}}^{\lambda(e_{-})}$

(5) $r = max^{\epsilon}H_{r}$ if and only if $\sigma_{0} = $$^{\lambda_{*}(e_{-})}w_{0}\phantom{~\hspace{-1pt}}^{\lambda_{*}(e_{+})}$

(6) $r = max^{\epsilon}L_{r}$ if and only if $\sigma_{-} = 1$, $\sigma_{0} = $$^{\lambda_{*}(e_{-})}w_{0}\phantom{~\hspace{-1pt}}^{\lambda_{*}(e_{+})}$

(7) $r = max^{\epsilon}R_{r}$ if and only if $\sigma_{0} = $$^{\lambda_{*}(e_{-})}w_{0}\phantom{~\hspace{-1pt}}^{\lambda_{*}(e_{+})}$, $\sigma_{+} = 1$

(8) $r = max^{\epsilon}J_{r}$ if and only if $\sigma_{-} = 1$, $\sigma_{0} = $$^{\lambda_{*}(e_{-})}w_{0}\phantom{~\hspace{-1pt}}^{\lambda_{*}(e_{+})}$, $\sigma_{+} = 1$
\end{thm}

\begin{proof}
The proofs are similar, so we will just look at $\varepsilon = +$.

(1) and (5) By definition of vanilla form we know $\sigma_{0} \in W_{\lambda^{*}(e_{-})}\Big($$^{\lambda_{*}(e_{-})}w_{0}^{\lambda_{*}(e_{+})}\Big)W_{\lambda^{*}(e_{+})} \subseteq W_{\lambda(e_{-})}w_{0}W_{\lambda(e_{+})}$. Then we can see $^{\lambda(e_{-})}w_{0}\phantom{~\hspace{-1pt}}^{\lambda(e_{+})} \leq \sigma_{0} \leq w_{0}$ and then, $^{\lambda(e_{-})}w_{0}\phantom{~\hspace{-1pt}}^{\lambda(e_{+})} \leq \sigma_{0} \leq~ $$^{\lambda_{*}(e_{-})}w_{0}\phantom{~\hspace{-1pt}}^{\lambda_{*}(e_{+})}$. By Corollary \ref{Green's relations and the Adherence order} we can then conclude that $\sigma_{-}e_{-}$$^{\lambda(e_{-})}w_{0}\phantom{~\hspace{-1pt}}^{\lambda(e_{+})}e_{+}\sigma_{+} \leq^{+} \sigma_{-}e_{-}\sigma_{0}e_{+}\sigma_{+} \leq^{+} \sigma_{-}e_{-}$$^{\lambda_{*}(e_{-})}w_{0}\phantom{~\hspace{-1pt}}^{\lambda_{*}(e_{+})}e_{+}\sigma_{+}$, each of which is written in vanilla form, and (by Corollary \ref{Green's relations and vanilla form}) belong to the same $\Hc$-class. As $r$ was arbitrary, we are done.

(2) and (6) Notice that for any $\sigma_{-}$ and $\sigma_{0}$ satisfying the conditions of vanilla form, $\sigma_{-}\sigma_{0} \in W^{\lambda_{*}(e_{+})} \subseteq W$. We can see that $1 \leq \sigma_{-}\sigma_{0} \leq w_{0}$ and thus, $1 \leq \sigma_{-}\sigma_{0} \leq w_{0}\phantom{~\hspace{-1pt}}^{\lambda_{*}(e_{+})} = $$^{\lambda_{*}(e_{-})}w_{0}\phantom{~\hspace{-1pt}}^{\lambda_{*}(e_{+})}$. It follows that $1e_{+}\sigma_{+} \leq^{+} (\sigma_{-}\sigma_{+})e_{+}\sigma_{+} \leq^{+} $$^{\lambda_{*}(e_{-})}w_{0}\phantom{~\hspace{-1pt}}^{\lambda_{*}(e_{+})}e_{+}\sigma_{+}$

If $r = min^{+}L_{r}$ then it is not hard to see $r = min^{+}H_{r}$. So by (1) we see that $\sigma_{0} = $$^{\lambda(e_{-})}w_{0}\phantom{~\hspace{-1pt}}^{\lambda(e_{+})}$, and hence $\sigma_{-} = \sigma_{+}^{-1} = $$^{\lambda(e_{+})}w_{0}\phantom{~\hspace{-1pt}}^{\lambda(e_{-})}$. Thus, $r = $$^{\lambda(e_{+})}w_{0}\phantom{~\hspace{-1pt}}^{\lambda(e_{-})}e_{-}$$^{\lambda(e_{-})}w_{0}\phantom{~\hspace{-1pt}}^{\lambda(e_{+})}e_{+}\sigma_{+}$. If $r = max^{+}L_{r}$ then it is not hard to see $r = max^{+}H_{r}$. So by (5) we see that $\sigma_{0} = $$^{\lambda_{*}(e_{-})}w_{0}\phantom{~\hspace{-1pt}}^{\lambda_{*}(e_{+})}$, and hence $\sigma_{-} = 1$. So $r = 1e_{-}$$^{\lambda_{*}(e_{-})}w_{0}\phantom{~\hspace{-1pt}}^{\lambda_{*}(e_{+})}e_{+}\sigma_{+}$.

$^{\lambda_{*}(e_{-})}w_{0}\phantom{~\hspace{-1pt}}^{\lambda_{*}(e_{+})} = w_{0}\phantom{~\hspace{-1pt}}^{\lambda_{*}(e_{+})}$ This demonstrates (6).

(3) and (7) are done similarly to (2) and (6).

(4) and (8) If $r = min^{+}J_{r}$ then it is clear that $r = min^{+}L_{r} = min^{+}R_{r} = min^{+}H_{r}$ and by combining our previous results, $r = ^{\lambda(e_{+})}w_{0}\phantom{~\hspace{-1pt}}^{\lambda(e_{-})}e_{-}$$^{\lambda(e_{-})}w_{0}\phantom{~\hspace{-1pt}}^{\lambda(e_{+})}e_{+}$$^{\lambda(e_{+})}w_{0}\phantom{~\hspace{-1pt}}^{\lambda(e_{-})}$.  If $r = max^{+}J_{r}$ then it is clear that $r = max^{+}L_{r} = max^{+}R_{r} = max^{+}H_{r}$ and by combining our previous results, $r = 1e_{-}$$^{\lambda_{*}(e_{-})}w_{0}\phantom{~\hspace{-1pt}}^{\lambda_{*}(e_{+})}e_{+}1$.
\end{proof}

In fact, we have already seen (and even without the finiteness assumption!) that $min^{+}L_{r}$ and $min^{-}R_{r}$ always exist. Indeed, the collection of such elements are $\GJ^{+}$ and $\JG^{-}$ respectively, from Theorem \ref{certain abs mins exist}. The elements introduced in Section 5 of \cite{PPR} are actually the minimum elements of $\Hc$-classes and hence form the familiar object, $\Oh$.

\begin{cor}
Let $\varepsilon = +$ or $-$. Define the following subsets of our monoid, $\Ren$: $\GJ^{\varepsilon} = \{r \in \Ren \mid r = min^{\varepsilon}L_{r}\}$, $\JG^{\varepsilon} = \{r \in \Ren \mid r = min^{\varepsilon}R_{r}\}$, $\NU^{\varepsilon} = \{r \in \Ren \mid r = min^{\varepsilon}J_{r}\}$, and
$\Oh^{\varepsilon} = \{r \in \Ren \mid r = min^{\varepsilon}H_{r}\}$.

(1) $\Oh^{+} = \Oh^{-}$ and we will henceforth refer to it as $\Oh$

(2) $\NU^{\varepsilon} = \GJ^{\varepsilon}\cap\JG^{\varepsilon}$

(3) each of these sets are monoids

(4) $\GJ^{+} = (\JG^{-})^{*}$ and $\JG^{+} =(\GJ^{-})^{*}$

(5) $\Oh$ is the smallest monoid containing $\GJ^{\varepsilon}$ and $\JG^{\varepsilon}$

(6) $\Oh$ is the smallest inverse monoid containing $\GJ^{\varepsilon}$

(7) $\Oh$ is the smallest inverse monoid containing $\JG^{\varepsilon}$

(8) $w_{0}\GJ^{\varepsilon} = \{r \in \Ren \mid r = max^{\varepsilon}L_{r}\}$

(9) $\JG^{\varepsilon}w_{0} = \{r \in \Ren \mid r = max^{\varepsilon}R_{r}\}$

(10) $w_{0}\Lambda = \Lambda^{-}w_{0} = \{r \in \Ren \mid r = max^{+}J_{r}\}$ and $w_{0}\Lambda^{-} = \Lambda w_{0} = \{r \in \Ren \mid r = max^{-}J_{r}\}$

(11) $w_{0}\Oh = \Oh w_{0} = \{r \in \Ren \mid r = max^{\varepsilon}H_{r}\}$ for each $\varepsilon$

(12) $\NU^{\varepsilon} \cong \Ren / \Jc$, $\GJ^{\varepsilon} \cong \Ren / \Lc$, $\JG^{\varepsilon} \cong \Ren / \Rc$, $\Oh \cong \Ren / \Hc$
\end{cor}

\begin{proof}
(1) This has already been covered by Theorem \ref{certain abs mins exist}, but is included for completeness.

(2) This is an application of Theorem \ref{abs elts vanilla conditions} (2), (3), and (4). Just comparing the vanilla forms in $(\Ren, \Lambda, S)$ gives the $\varepsilon = +$ case and a symmetry argument takes care of the $\varepsilon = -$ case.

(3) Each of $\Oh$, $\GJ^{+}$, and $\JG^{-}$ are monoids by Theorem \ref{certain abs mins exist}. The monoidness of $\JG^{+}$ and $\GJ^{-}$ are simiarly demonstrated, but as part of the Renner-Coxeter system, $(\Ren, \Lambda^{-}, S)$. By (2), $\NU^{+}$ and $\NU^{-}$ are intersections of monoids, and hence monoids themselves.

(4) The first result is Theorem \ref{certain abs mins exist} (2) and the second follows in $(\Ren, \Lambda^{-}, S)$. 

(5) By Theorem \ref{certain abs mins exist} we know $\Oh$ is the smallest monoid containing $\GJ^{+}$ and $\JG^{-}$. By similar reasoning, but in the setting of the $(\Ren, \Lambda^{-}, S)$ system we can conclude that $\Oh$ is also the smallest monoid containing $\GJ^{-}$ and $\JG^{+}$. Combining these results gives us (5).

(6) For $\varepsilon = +$ this has been done in Theorem \ref{certain abs mins exist} (4), and the case for $-$ follows by similar reasoning, but in $(\Ren, \Lambda^{-}, S)$.

(7) For $\varepsilon = -$ this has been done in Theorem \ref{certain abs mins exist} (4), and the case for $+$ follows by similar reasoning, but in $(\Ren, \Lambda^{-}, S)$.

(8) Consider $r = w_{0}ey$ for any $ey \in \GJ^{+}$. Then $r = w_{0}^{\lambda_{*}(e)}ey$ is written in left standard form. For any $s \Lc r$, we can write $s = xey$ for $x \in W^{\lambda_{*}(e)}$. But, $x \in W^{\lambda_{*}(e)}$ implies (by Theorem \ref{Coxeter -fixes}) that $x \leq_{L} w_{0}^{\lambda_{*}(e)}$ and so $x \leq w_{0}^{\lambda_{*}(e)}$. But this implies that $xey \leq^{+} w_{0}^{\lambda_{*}(e)}ey$ and so $s \leq^{+} r$. Since $s$ was arbitrary, we see $r = max^{+}L_{r}$, and since every $\Lc$-class in $\Ren$ has exactly one element of $\GJ^{+}$ in it (by definition) then every $\Lc$-class has such a maximum element, concluding our result.

The $\varepsilon = -$ case is handled analogously.  

(9) is taken care of in a manner similar to (8).

(10) For any $\Jc$-class, we see that there is exactly one element of the form $w_{0}\Lambda$, so it suffices to show that this element is largest in the given $\Jc$-class. Let $w_{0}e = w_{0}^{\lambda_{*}(e)}e$ be such an element written in left standard form. Consider an arbitrary element $xey$ written in left standard form from the same $\Jc$-class. Then $x \in W^{\lambda_{*}(e)}$ implies that $x \leq w_{0}^{\lambda_{*}(e)}$. Adding to this the fact that $1 \leq y$, we conclude that $xey \leq^{+} w_{0}e$. So $w_{0}\Lambda = \{r \in \Ren \mid r = max^{+}J_{r}\}$. Simiarly we can show that $\Lambda^{-}w_{0} = \{r \in \Ren \mid r = max^{+}J_{r}\}$.

The $\varepsilon = -$ case is completed by performing the same actions, but in $(\Ren, \Lambda^{-}, S)$, allowing us to conclude that (with regards to $(\Ren, \Lambda, S)$) $w_{0}\Lambda^{-} = \Lambda w_{0} =  \{r \in \Ren \mid r = max^{-}J_{r}\}$.

(11) Consider any $r \in \Oh$ and $s \in H_{r}$. Then $r \leq^{\varepsilon} s$. By Corollary \ref{Adherence equivalences} we know that $w_{0}rw_{0} \leq^{-\varepsilon} w_{0}sw_{0}$. Since $s \in H_{r}$ $\Leftrightarrow$ $w_{0}sw_{0} \in H_{w_{0}rw_{0}}$ we conclude that $w_{0}rw_{0} \in \Oh$ as well. But then $w_{0}\Oh w_{0} = \Oh$ and hence $w_{0}\Oh = \Oh w_{0}$ It remains to show that these elements are maximum in their $\Hc$-classes.

Consider an arbitrary $r \in \Oh$ and any element $s \in H_{r}$. Let $r = xey$ and $s = zey$ be written in left standard form. Then $w_{0}r = (w_{0}x)^{\lambda^{*}(e)}ey$ and $w_{0}s = (w_{0}z)^{\lambda_{*}(e)}ey$, So it is clear that $w_{0}s \Hc w_{0}r$ if and only if $s \Hc r$. We can observe the following equivalences, $r \leq^{+} s$ if and only if $x \leq z$ if and only if $w_{0}z \leq w_{0}x$ if and only if $(w_{0}z)^{\lambda_{*}(e)} \leq (w_{0}x)^{\lambda_{*}(e)}$ if and only if $w_{0}s \leq^{+} w_{0}r$. This allows us to conclude that $w_{0}\Oh = \{r \in \Ren \mid r = max^{+}H_{r}\}$. A similar chain of equivalences shows $\Oh w_{0} = \{r \in \Ren \mid r = max^{-}H_{r}\}$.

(12) This is similar to Theorem \ref{certain abs mins exist}. We will just handle the $\NU^{+}$ case, as the rest have similar reasoning. From the definition it is clear that $r, s \in \NU^{+}$ and $r \Jc s$ implies that $r \leq^{+} s$ and $s \leq^{+} r$, hence $r = s$. It remains to show that each $\Jc$-class has a minimum. However, in the preceeding theorem, we contructed $min^{+}J_{r}$ by its vanilla form from an arbitrary element $r \in \Ren$. So since for all $r \in \Ren$ we have $\abs{J_{r} \cap \NU^{+}} = 1$ we conclude $\NU^{+} \cong \Ren / \Jc$.
\end{proof}

In fact, the proof of our last result can be extended to show the following corollary.

\begin{cor}
Let $\varepsilon = +$ or $-$ and let $r \leq^{\varepsilon} s \in \Ren$

(1) If $r \Lc s$ then $w_{0}r \Lc w_{0}s$ and $w_{0}s \leq^{\varepsilon} w_{0}r$

(2) If $r \Rc s$ then $rw_{0} \Rc sw_{0}$ and $sw_{0} \leq^{\varepsilon} rw_{0}$
\end{cor}

\begin{proof}
We shall just prove (1) with $\varepsilon = +$ as the others are done simiarly. Let $r = xey$ and $s = zey$ be written in left standard form. Then $w_{0}r = (w_{0}x)^{\lambda^{*}(e)}ey$ and $w_{0}s = (w_{0}z)^{\lambda_{*}(e)}ey$, So it is clear that $w_{0}s \Lc w_{0}r$ since $r \Lc s$. We can observe the following equivalences, $r \leq^{+} s$ if and only if $x \leq z$ if and only if $w_{0}z \leq w_{0}x$ if and only if $(w_{0}z)^{\lambda_{*}(e)} \leq (w_{0}x)^{\lambda_{*}(e)}$ if and only if $w_{0}s \leq^{+} w_{0}r$. This concludes the result.
\end{proof}

The absolute elements of a $\Tc$-class can in some sense be considered as indicators of the ordering of the $\Tc$-classes in the whole poset, $(\Ren, \leq^{\varepsilon})$. This next theorem showcases this result.

\begin{thm} \label{absolute equivalences}
Let $r, s \in \Ren$ be arbitrary elements and suppose $\Tc$ is an equivalence relation. Consider the following statements:

(1) $min^{\varepsilon}T_{r} \leq^{\varepsilon} min^{\varepsilon}T_{s}$

(2) $max^{\varepsilon}T_{r} \leq^{\varepsilon} max^{\varepsilon}T_{s}$

(3) There exist $a \in T_{r}$ and $b \in T_{s}$ so that $a \leq^{\varepsilon} b$

\noindent If $\Tc = \Jc$, $\Lc$, or $\Rc$ then all three are equivalent. If $\Tc = \Hc$ then only (2) and (3) are equivalent.
\end{thm}

\begin{proof}
Since $min^{\varepsilon}T_{r} \Tc r \Tc max^{\varepsilon}T_{r}$ for all $r \in \Ren$ we can quickly see that (2) implies (3).

Case $\Tc = \Jc$: Let $r = \sigma_{-}e_{-}\sigma_{0}e_{+}\sigma_{+}$ and $s = \tau_{-}f_{-}\tau_{0}f_{+}\tau_{+}$ be written in vanilla form. By using Theorem \ref{abs elts vanilla conditions} we can notice $min^{+}J_{r} = \Big($$^{\lambda(e_{+})}w_{0}^{\lambda(e_{-})}\Big)e_{-}\Big($$^{\lambda(e_{-})}w_{0}^{\lambda(e_{+})}\Big)e_{+}\Big($$^{\lambda(e_{+})}w_{0}^{\lambda(e_{-})}\Big)$, $min^{+}J_{s} = \Big($$^{\lambda(f_{+})}w_{0}^{\lambda(f_{-})}\Big)f_{-}\Big($$^{\lambda(f_{-})}w_{0}^{\lambda(f_{+})}\Big)f_{+}\Big($$^{\lambda(f_{+})}w_{0}^{\lambda(f_{-})}\Big)$, $max^{+}J_{r} = 1e_{-}\Big($$^{\lambda_{*}(e_{-})}w_{0}^{\lambda_{*}(e_{+})}\Big)e_{+}1$, and $max^{+}J_{s} = 1f_{-}\Big($$^{\lambda_{*}(f_{-})}w_{0}^{\lambda_{*}(f_{+})}\Big)f_{+}1$.

Suppose (1) is true. Then $min^{+}J_{r} \leq^{+} min^{+}J_{s}$ implies that $e_{-} \leq f_{-}$ and $e_{+} \leq f_{+}$. It follows that $\lambda_{*}(f_{-}) \subseteq \lambda_{*}(e_{-})$ and $\lambda_{*}(f_{+}) \subseteq \lambda_{*}(e_{+})$ and so, $\Big($$^{\lambda_{*}(e_{-})}w_{0}^{\lambda_{*}(e_{+})}\Big) \leq \Big($$^{\lambda_{*}(f_{-})}w_{0}^{\lambda_{*}(f_{+})}\Big)$. But it is exactly this condition that shows us that $max^{+}J_{r} \leq^{+} max^{+}J_{s}$.

Suppose (3). Without loss of generality we may assume $a = r$ and $b = s$. It follows from $r \leq^{+} s$ that $e_{-} \leq f_{-}$ (hence $\lambda_{*}(f_{-}) \subseteq \lambda_{*}(e_{-})$ and $\lambda^{*}(e_{-}) \subseteq \lambda^{*}(f_{-})$) and $e_{+} \leq f_{+}$ (hence $\lambda_{*}(f_{+}) \subseteq \lambda_{*}(e_{+})$ and $\lambda^{*}(e_{+}) \subseteq \lambda^{*}(f_{+})$). By Theorem \ref{Coxeter -fixes} we can see that there exist $c \in W_{\lambda^{*}(f_{+})}$, $d \in W_{\lambda_{*}(f_{+})}$ so that $^{\lambda(f_{-})}w_{0}^{\lambda(f_{+})} = w_{0}^{\lambda(f_{+})} = w_{0}cd$, exist $a \in W_{\lambda_{*}(e_{+})}$, $b \in W_{\lambda^{*}(e_{+})}$ so that  $^{\lambda(e_{-})}w_{0}^{\lambda(e_{+})} = w_{0}^{\lambda(e_{+})} = w_{0}ab$, exist $\gamma \in W_{\lambda^{*}(f_{-})}$, $\delta \in W_{\lambda_{*}(f_{-})}$ so that $^{\lambda(f_{-})}w_{0}^{\lambda(f_{+})} = $$^{\lambda(f_{-})}w_{0} = \delta\gamma w_{0}$, and exist $\alpha \in W_{\lambda_{*}(e_{-})}$, $\beta \in W_{\lambda^{*}(e_{-})}$ so that $^{\lambda(e_{-})}w_{0}^{\lambda(e_{+})} = $$^{\lambda(e_{-})}w_{0} = \beta\alpha w_{0}$.

Then we can observe (after some minor rearrangement of the above) that $^{\lambda(f_{+})}w_{0}^{\lambda(f_{-})}(\gamma^{-1}\beta)(\delta^{-1}\alpha) = $$^{\lambda(e_{+})}w_{0}^{\lambda(e_{-})}$, $^{\lambda(e_{-})}w_{0}^{\lambda(e_{+})} = (\alpha{-1}\delta)(\beta^{-1}\gamma)\Big(^{\lambda(f_{-})}w_{0}^{\lambda(f_{+})}\Big)(cb^{-1})(da^{-1})$, and $(ad^{-1})(bc^{-1})$$^{\lambda(f_{+})}w_{0}^{\lambda(f_{-})} = $$^{\lambda(e_{+})}w_{0}^{\lambda(e_{-})}$, which tells us that $min^{+}J_{r} \leq min^{+}J_{s}$.

Case $\Tc = \Lc$: Let $r = \sigma e_{+}\sigma_{+}$ and $s = \tau f_{+}\tau_{+}$ be written in left standard form. By using Theorem \ref{abs elts vanilla conditions} we can notice $min^{+}L_{r} = 1 e_{+}\sigma_{+}$, $min^{+}L_{s} = 1f_{+}\tau_{+}$, $max^{+}L_{r} = w_{0}^{\lambda_{*}(e_{+})}e_{+}\sigma_{+}$, and $max^{+}L_{s} = w_{0}^{\lambda_{*}(f_{+})}f_{+}\tau_{+}$.

Suppose (1) is true. Then $min^{+}L_{r} \leq^{+} min^{+}L_{s}$ implies that $e_{+} \leq f_{+}$ and there exists $w_{+} \in W_{\lambda^{*}(f_{+})}W_{\lambda_{*}(e_{+})}$ so that $1 \leq 1w_{+}$ and $w_{+}^{-1}\tau_{+} \leq \sigma_{+}$. It follows that $\lambda_{*}(f_{+}) \subseteq \lambda_{*}(e_{+})$ and so, $w_{0}^{\lambda_{*}(e_{+})} \leq w_{0}^{\lambda_{*}(f_{+})}$. Let $w_{+} = uv$ where $u \in W_{\lambda_{*}(f_{+})}$ and $v \in W_{\lambda_{*}(e_{+})}$. Then $^{\lambda_{*}(e_{+})}(\tau_{+}) \leq $$^{\lambda_{*}(e_{+})}(u^{-1}\tau_{+}) \leq v^{-1}u^{-1}\tau_{+} \leq \sigma_{+}$ so we can find $w \in W_{\lambda_{*}(e_{+})}$ so that $w^{-1}\tau_{+} \leq \sigma_{+}$. Noting as well that $w_{0}^{\lambda_{*}(e_{+})} \leq w_{0}^{\lambda_{*}(f_{+})}w$ it follows that $max^{+}L_{r} \leq^{+} max^{+}L_{s}$.

Suppose (3). Without loss of generality we may assume $a = r$ and $b = s$. It follows from $r \leq^{+} s$ that $e_{+} \leq f_{+}$ and there exists $w_{+} \in W_{\lambda^{*}(f_{+})}W_{\lambda_{*}(e_{+})}$ so that $\sigma \leq \tau w_{+}$ and $w_{+}^{-1}\tau_{+} \leq \sigma_{+}$. It is a quick follow-up to observe that $1 \leq w_{+} = 1w_{+}$ and $w_{+}^{-1}\tau_{+} \leq \sigma_{+}$, thus $min^{+}L_{r} \leq min^{+}L_{s}$.

Case $\Tc = \Rc$: Is nearly identical to the $\Lc$ case.

Case $\Tc = \Hc$: Let $r = \sigma_{-}e_{-}\sigma_{0}e_{+}\sigma_{+}$ and $s = \tau_{-}f_{-}\tau_{0}f_{+}\tau_{+}$ be written in vanilla form. By using Theorem \ref{abs elts vanilla conditions} we can notice $min^{+}H_{r} = \sigma_{-}e_{-}\Big($$^{\lambda(e_{-})}w_{0}^{\lambda(e_{+})}\Big)e_{+}\sigma_{+}$, $min^{+}H_{s} = \tau_{-}f_{-}\Big($$^{\lambda(f_{-})}w_{0}^{\lambda(f_{+})}\Big)f_{+}\tau_{+}$, $max^{+}H_{r} = \sigma_{-}e_{-}\Big($$^{\lambda_{*}(e_{-})}w_{0}^{\lambda_{*}(e_{+})}\Big)e_{+}\sigma_{+}$, and $max^{+}H_{s} = \tau_{-}f_{-}\Big($$^{\lambda_{*}(f_{-})}w_{0}^{\lambda_{*}(f_{+})}\Big)f_{+}\tau_{+}$.

Suppose (3). As we have been doing, we will just assume $a = r$ and $b = s$. We can find $w_{-} \in W_{\lambda_{*}(e_{-})}W_{\lambda^{*}(f_{-})}$ and $w_{+} \in W_{\lambda^{*}(f_{+})}W_{\lambda_{*}(e_{+})}$ so that $\tau_{-}w_{-}^{-1} \leq \sigma_{-}$ and $w_{+}^{-1}\tau_{+} \leq \sigma_{+}$. Like above we can assume that $w_{-} \in W_{\lambda_{*}(e_{-})}$ and $w_{+} \in W_{\lambda_{*}(e_{+})}$. As well, we also know that $\lambda_{*}(f_{-}) \subseteq \lambda_{*}(e_{-})$ and $\lambda_{*}(f_{+}) \subseteq \lambda_{*}(e_{+})$. Thus $^{\lambda_{*}(e_{-})}w_{0}^{\lambda_{*}(e_{+})} \leq $$^{\lambda_{*}(f_{-})}w_{0}^{\lambda_{*}(f_{+})}$ and hence $^{\lambda_{*}(e_{-})}w_{0}^{\lambda_{*}(e_{+})} \leq w_{-}$$^{\lambda_{*}(f_{-})}w_{0}^{\lambda_{*}(f_{+})}w_{+}$. This shows us that $max^{+}H_{r} \leq^{+} max^{+}H_{s}$.
\end{proof}

The inspiration for our look at minimum elements in equivalence classes comes from Section 5 in \cite{PPR} where a number of results, noteably Corollary 5.5, deal with the properties of the elements of $\mathcal{O}$ inside their $\Hc$-classes. Were it not for this section, the author's PhD thesis would not have taken the direction it did, and consequently this paper would not have been written as a combinatorial follow up. 

It is somewhat ironic that the result this work was inspired by, Corollary 5.5 in \cite{PPR}, actually states that $r \leq^{+} s$ implies $min^{+}H_{r} \leq^{+} min^{+}H_{s}$, which one may notice was exactly the statement that had to be left out of the preceeding theorem, as it is not correct.

\begin{eg}
Consider the elements of the $3\times 3$ Rook monoid with the usual Adherence order given by the invertible upper triangular matrices, $r = $\scalebox{0.5}{$\left(\begin{array}{ccc}0 & 0 & 0 \\ 0 & 1 & 0 \\ 1 & 0 & 0\end{array}\right)$} and $s = $\scalebox{0.5}{$\left(\begin{array}{ccc}0 & 0 & 1 \\ 0 & 1 & 0 \\ 1 & 0 & 0\end{array}\right)$}. Then $r \leq^{+} s$ but $min^{+}H_{r} \not{\leq^{+}} min^{+}H_{s}$. Indeed, $min^{+}H_{r} = $\scalebox{0.5}{$\left(\begin{array}{ccc}0 & 0 & 0 \\ 1 & 0 & 0 \\ 0 & 1 & 0\end{array}\right)$} and $min^{+}H_{s} = $\scalebox{0.5}{$\left(\begin{array}{ccc}1 & 0 & 0 \\ 0 & 1 & 0 \\ 0 & 0 & 1\end{array}\right)$}.
\end{eg}

As a consequence of this theorem, we can show that chains in $\Ren$ do not wander in and out of equivalence classes. In a sense, equivalence classes are their own self-contained bubbles when it comes to our Adherence orders. 

\begin{cor}
Suppose that we have an equivalence relation $\Tc = \Jc$, $\Lc$, $\Rc$, or $\Hc$, and a chain of elements, $r_{0} \leq^{\varepsilon} r_{1} \leq^{\varepsilon} \cdots \leq^{\varepsilon} r_{n}$, in $\Ren$. Then $r_{0} \Tc r_{n}$ implies that $r_{0} \Tc r_{i}$ for all $i$.
\end{cor}

\begin{proof}
Using Theorem \ref{absolute equivalences} we can see that our chain condition implies $max^{\varepsilon}T_{r_{0}} \leq^{\varepsilon} max^{\varepsilon}T_{r_{1}} \leq^{\varepsilon} \cdots \leq^{\varepsilon} max^{\varepsilon}T_{r_{n}}$. But since $r_{0} \Tc r_{n}$ we see that $T_{r_{0}} = T_{r_{n}}$ and hence $max^{\varepsilon}T_{r_{0}} \leq^{\varepsilon} max^{\varepsilon}T_{r_{1}} \leq^{\varepsilon} \cdots \leq^{\varepsilon} max^{\varepsilon}T_{r_{0}}$. So by squeezing we see that $min^{\varepsilon}T_{r_{0}} = max^{\varepsilon}T_{r_{i}}$ for all $i$. By definition, $s \Tc max^{\varepsilon}T_{s}$ for all $s \in \Ren$ we conclude that $r_{0} \Tc r_{i}$ for all $i$.
\end{proof}

\section{Geometric Interpretations} \label{geometry}

The scope of the generalized Renner-Coxeter systems is indeed impressive and allows us to consider the theory of monoids of Lie type, Renner monoids of reductive algebraic monoids, and Renner monoids created from Kac-Moody groups. However, let us take a moment and look at the implications of our definitions and results when confined to Renner monoids of reductive algebraic monoids.

Such Renner monoids are always finite, so everything we have dealt with in this paper applies. In a Renner monoid, the Adherence order corresponds exactly to our definition of $\leq^{+}$ (Corollary 1.5 in \cite{PPR}). This fact, along with the knowledge that all Renner monoids coming from reductive monoids are finite (Proposition 3.2.1 in \cite{Renner Bruhat I}) means that the vanilla form and all the results we have derived from it now apply.

For the ease of our notation, we will rebrand our Borel subgroup as $B^{+} = B$ and our irreducible reductive monoid, $M$ with group of units $G$. 

\begin{thm} \label{geometric interpretation}
~

(1) $B^{\varepsilon}T_{r}B^{\varepsilon} = \bigsqcup_{s \in T_{r}}B^{\varepsilon}sB^{\varepsilon}$ is an irreducible variety. 

(2) $B^{\varepsilon}max^{\varepsilon}T_{r}B^{\varepsilon}$ is the unique dense orbit of $B^{\varepsilon}\times B^{\varepsilon}$ on $B^{\varepsilon}T_{r}B^{\varepsilon}$.

(3) $B^{\varepsilon}min^{\varepsilon}T_{r}B^{\varepsilon}$ is the unique closed orbit of $B^{\varepsilon}\times B^{\varepsilon}$ on $B^{\varepsilon}T_{r}B^{\varepsilon}$.
\end{thm}

\begin{proof}
(1) The individual cases when $\varepsilon = +$ are covered in \cite{O'Hara Thesis} in Section 4. The $B^{-}$ cases can be derived through nearly identical reasoning.

(2) By definition we know, $max^{\varepsilon}T_{r} \Tc r$ so we can see that $B^{\varepsilon}max^{\varepsilon}T_{r}B^{\varepsilon} \subseteq B^{\varepsilon}T_{r}B^{\varepsilon}$. It follows from there, $\overline{B^{\varepsilon}max^{\varepsilon}T_{r}B^{\varepsilon}} \subseteq \overline{B^{\varepsilon}T_{r}B^{\varepsilon}} = \bigcup_{s \in T_{r}}\overline{B^{\varepsilon}sB^{\varepsilon}} \subseteq \bigcup_{s \in T_{r}}\overline{B^{\varepsilon}max^{\varepsilon}T_{r}B^{\varepsilon}} = \overline{B^{\varepsilon}max^{\varepsilon}T_{r}B^{\varepsilon}}$.

(3) Suppose that $s \in \overline{B^{\varepsilon}min^{\varepsilon}T_{r}B^{\varepsilon}}\cap B^{\varepsilon}T_{r}B^{\varepsilon} = \bigsqcup_{t \in T_{r}}\overline{B^{\varepsilon}min^{\varepsilon}T_{r}B^{\varepsilon}}\cap B^{\varepsilon}tB^{\varepsilon}$. If $s \in \overline{B^{\varepsilon}min^{\varepsilon}T_{r}B^{\varepsilon}}\cap B^{\varepsilon}tB^{\varepsilon}$ then $\overline{B^{\varepsilon}min^{\varepsilon}T_{r}B^{\varepsilon}}\cap B^{\varepsilon}tB^{\varepsilon} \neq \emptyset$, and hence $t \leq^{\varepsilon} min^{\varepsilon}T_{r}$. However, by definition this means $t = min^{\varepsilon}T_{r}$. Thus $\overline{B^{\varepsilon}min^{\varepsilon}T_{r}B^{\varepsilon}}\cap B^{\varepsilon}T_{r}B^{\varepsilon} = B^{\varepsilon}min^{\varepsilon}T_{r}B^{\varepsilon}$, hence $B^{\varepsilon}min^{\varepsilon}T_{r}B^{\varepsilon}$ is a closed orbit. By similar analysis, we can see that for any other orbit in $B^{\varepsilon}T_{r}B^{\varepsilon}$, $B^{\varepsilon}min^{\varepsilon}T_{r}B^{\varepsilon}$ will always be contained in its closure, making $B^{\varepsilon}min^{\varepsilon}T_{r}B^{\varepsilon}$ the only closed orbit.
\end{proof}

\begin{prop} \label{extension of <T conditions}
Let us define a new relation on the $\Tc$-classes of $\Ren$. $T_{r} \leq^{\varepsilon}_{\Tc} T_{s}$ if and only if $B^{\varepsilon}T_{r}B^{\varepsilon} \subseteq \overline{B^{\varepsilon}T_{r}B^{\varepsilon}}$. Like $B^{\varepsilon}rB^{\varepsilon} \subseteq \overline{B^{\varepsilon}sB^{\varepsilon}}$, $\leq^{\varepsilon}_{\Tc}$ can be shown to be a partial order on the $\Tc$-classes of $\Ren$. Consider the following statements:

(1) $min^{\varepsilon}T_{r} \leq^{\varepsilon} min^{\varepsilon}T_{s}$

(2) $max^{\varepsilon}T_{r} \leq^{\varepsilon} max^{\varepsilon}T_{s}$

(3) there exist $a \in T_{r}$ and $b \in T_{s}$ with $a \leq^{\varepsilon} b$

(4) $T_{r} \leq^{\varepsilon}_{\Tc} T_{s}$

\noindent If $\Tc = \Jc$, $\Lc$, or $\Rc$ then all four are equivalent. If $\Tc = \Hc$ the only the last three are equivalent.
\end{prop}

\begin{proof}
The equivalence of (2) and (3) (and (1) is $\Tc \neq \Hc$) has already been established by Theorem \ref{absolute equivalences}. Suppose (2). Then $max^{\varepsilon}T_{r} \leq^{\varepsilon} max^{\varepsilon}T_{s}$ if and only if $B^{\varepsilon}max^{\varepsilon}T_{r}B^{\varepsilon} \subseteq \overline{B^{\varepsilon}max^{\varepsilon}T_{s}B^{\varepsilon}}$. However, by Theorem \ref{geometric interpretation} we can then see that this implies $B^{\varepsilon}T_{r}B^{\varepsilon} \subseteq \overline{B^{\varepsilon}T_{r}B^{\varepsilon}} = \overline{B^{\varepsilon}max^{\varepsilon}T_{r}B^{\varepsilon}} \subseteq \overline{B^{\varepsilon}max^{\varepsilon}T_{s}B^{\varepsilon}} = \overline{B^{\varepsilon}T_{s}B^{\varepsilon}}$ which is (4).

Supposing (4) then, $\bigsqcup_{u \Tc r} B^{\varepsilon}uB^{\varepsilon} = B^{\varepsilon}T_{r}B^{\varepsilon} \subseteq \overline{B^{\varepsilon}T_{s}B^{\varepsilon}} = \overline{\bigsqcup_{v\Tc s}B^{\varepsilon}vB^{\varepsilon}} = \bigsqcup_{v\Tc s}\overline{B^{\varepsilon}vB^{\varepsilon}} = \bigsqcup_{w\leq^{\varepsilon} v, v\Tc s}B^{\varepsilon}wB^{\varepsilon}$. This means we can find $a \Tc r$ and $c \leq^{\varepsilon} b \Tc s$ so that $BaB \cap BcB \neq \emptyset$. However, this condition is true exactly when $a = c$ (Theorem 8.3 in \cite{Renner Book}) and so $a \leq^{\varepsilon} b$, showing (3).
\end{proof}

$\leq^{\varepsilon}_{\Tc}$ represents a generalization of $\leq^{\varepsilon}$. Indeed, for the equivalence relation of equality, we see $\leq^{\varepsilon} = \leq^{\varepsilon}_{=}$. This Adherence order on $\Tc$-classes lets us bridge the gap between equality and the usual partial order on $\Jc$-classes, as both are examples of these $\leq^{\varepsilon}_{\Tc}$ relations. Indeed, $J_{r} \subseteq \Ren J_{s}\Ren$ if and only if $J_{r} \leq^{+}_{\Jc} J_{s}$ if and only if $J_{r} \leq^{-}_{\Jc} J_{s}$.

Using our Borel subgroups, we can recognize our absolute minimum elements in a more geometric way.

\begin{prop} \label{geometric definitions for min elts}
Pick $\varepsilon = +$ or $-$. Then,

(1) $\GJ^{\varepsilon} = \{r \in \Ren \mid B^{\varepsilon}r \subseteq rB^{\varepsilon}\} = \{r \in \Ren \mid B^{\varepsilon}rB^{\varepsilon} = rB^{\varepsilon}\}$

(2) $\JG^{\varepsilon} = \{r \in \Ren \mid rB^{\varepsilon} \subseteq B^{\varepsilon}r\} = \{r \in \Ren \mid B^{\varepsilon}rB^{\varepsilon} = B^{\varepsilon}r\}$

(3) $\NU^{\varepsilon} = \{r \in \Ren \mid B^{\varepsilon}r = rB^{\varepsilon}\}$

(4) $\Oh^{\varepsilon} = \{r \in \Ren \mid r^{*}B^{\varepsilon}r \subseteq r^{*}rB^{\varepsilon}\} = \{r \in \Ren \mid rB^{\varepsilon}r^{*} \subseteq B^{\varepsilon}rr^{*}\}$
\end{prop}

\begin{proof}
(1) and (2) $\GJ^{+}$ can be found to be equivalent to both sets by Section 9 in \cite{Renner Bruhat I}. The other three sets of absolute minimums are found by analogous work. 

(3) This is a combination of (1) and (2) as $\NU^{\varepsilon} = \GJ^{\varepsilon}\cap\JG^{\varepsilon}$.

(4) It is clear from the definitions of these sets that $\GJ^{+}$, $\GJ^{-}$, $\JG^{+}$, and $\JG^{-}$ are subsets. It is also clear that these two sets are monoids. In \cite{Renner Bruhat II} with Corollary 2.3 these sets are shown to be equal and further results show that they form a monoid. By Theorem 2.8 of that same paper, it is noted that these sets are equal to $\Ren / \Hc$. However, the $\Oh$ we defined in this paper is the smallest monoid containing $\GJ^{+}$, $\GJ^{-}$, $\JG^{+}$, and $\JG^{-}$ and also has the property that $\Oh = \Ren / \Hc$. It follows that these sets indeed describe $\Oh$.
\end{proof}

Infact, it is Renner's work in \cite{Renner Bruhat I} from which $\GJ^{+}$ gets is symbol, as these are the {\bf Gauss-Jordan} elements of the Renner monoid. The symbol $\GJ$ appears in other work, such as \cite{Renner Descent} where it is used to define $\JG^{+}$ instead. This is only a cosmetic difference, but is still noteworthy to those who wish to read both papers. Simiarly, $\Oh$ gets its symbol from \cite{Renner Bruhat II} as it contains the set of {\bf order-preserving} elements of $\Ren$.

Combining this most recent result with Proposition \ref{extension of <T conditions} we can simplify the assessment of $\leq_{\Lc}$ and $\leq_{\Rc}$.

\begin{cor}
~

(1) For $r, s \in \GJ^{\varepsilon}$ then $L_{r} \leq^{\varepsilon}_{\Lc} L_{s}$ if and only if $rB^{\varepsilon} \subseteq \overline{sB^{\varepsilon}}$.

(2) For $r, s \in \JG^{\varepsilon}$ then $R_{r} \leq^{\varepsilon}_{\Rc} R_{s}$ if and only if $B^{\varepsilon}r \subseteq \overline{B^{\varepsilon}s}$.
\end{cor}

\begin{proof}
(1) $L_{r} \leq^{\varepsilon}_{\Lc} L_{s}$ if and only if $r = min^{\varepsilon}L_{r} \leq^{\varepsilon} min^{\varepsilon}L_{s} = s$ if and only if $rB^{\varepsilon} = B^{\varepsilon}rB^{\varepsilon} \subseteq \overline{B^{\varepsilon}sB^{\varepsilon}} = \overline{sB^{\varepsilon}}$. (2) is done similarly.
\end{proof}

There is certainly much more work that can be done with Renner-Coxeter systems. They provide a solely combinatorial platform from which to investigate a traditionally algebraic and geometric structure. As a generalization of Coxeter systems, the overlap and extensions of results from \cite{Bjorner and Brenti} are definitely worth exploring and cataloguing.

The author wishes to thank Lex Renner for his time in discussions and his guidance on the creation of this paper and the research that lead to its completion.

\end{document}